\documentclass[12pt, dvipsnames, reqno]{amsart}
\usepackage[pagewise]{lineno}
\usepackage[T1]{fontenc}
\usepackage{a4wide}
\usepackage{dsfont}
\usepackage{yhmath,mathrsfs,amsthm,amsmath,amssymb,amsfonts, enumerate,lipsum, appendix,mathtools, bm}
\usepackage{bbold}
\usepackage[mathscr]{euscript}
\usepackage[normalem]{ulem} 
\usepackage{graphicx}
\usepackage{pgf,tikz}
\usepackage{tkz-euclide}
\usepackage{graphicx}
\usepackage{subcaption}
\usetikzlibrary{shapes.geometric}
\usetikzlibrary{arrows,hobby}
\usepackage{enumitem}
\usepackage[english]{babel}
\DeclareMathOperator\supp{supp}

\allowdisplaybreaks

\usepackage[sort,numbers]{natbib}
\usepackage[colorlinks=true,urlcolor=mygreen, citecolor=mygreen,linkcolor=blue,
linktocpage,pdfpagelabels, bookmarksnumbered,bookmarksopen]{hyperref}
\definecolor{mygreen}{rgb}{0.01,0.6,0.2}
\usepackage[margin=0.80in]{geometry}

\allowdisplaybreaks
\usepackage{orcidlink}
\usepackage{bigints}
 \usepackage{esint}
\newtheorem{theorem}{Theorem}[section]
\newtheorem{lemma}[theorem]{Lemma}
\newtheorem{proposition}[theorem]{Proposition}

\newtheorem{definition}[theorem]{Definition}
\theoremstyle{definition}

\newtheorem{remark}[theorem]{Remark}
\numberwithin{equation}{section}
\usepackage[hyperpageref]{backref}

\newcommand*\rd{\mathbb{R}^d}

\newcommand*\N{\mathcal{N}}
\newcommand{\al} {\alpha}
\newcommand{\pa} {\partial}
\newcommand{\be} {\beta}

\newcommand{\la} {\lambda}

\newcommand{\Gr} {\nabla}
\newcommand{\no} {\nonumber}
\newcommand{\noi} {\noindent}
\newcommand{\var} {\varepsilon}
\newcommand{\ra} {\rightarrow}

\newcommand{\bee} {\begin{equation}}
	\newcommand{\eee} {\end{equation}}
\newcommand{\bea} {\begin{eqnarray}}
	\newcommand{\eea} {\end{eqnarray}}
\newcommand{\Bea} {\begin{eqnarray*}}
	\newcommand{\Eea} {\end{eqnarray*}}
\def\Rn{\mathbb{R}^d}

\def\d{\,{\rm d}}
\def\dx{\,{\rm d}x}
\def\dy{\,{\rm d}y}

\def\C{{\mathcal C}}
\def\D{{\mathcal D}}

\def\W{{\mathcal W}}

\def\R{{\mathbb R}}
\def\N{{\mathbb N}}

\def\({{\Big(}}
\def\){{\Big)}}
\def\cc{{\C_c^\infty}}

\def\dx{\,{\rm d}x}
\def\dxy{\,{\rm d}x{\rm d}y}

\def\dh{\,{\rm d}h}
\DeclarePairedDelimiter\abs{\lvert}{\rvert}%
\DeclarePairedDelimiter\norm{\lVert}{\rVert}%
\def\wps{{\mathcal{D}^{s,p}}}
\usepackage{orcidlink}

\title[$p$-Fractional weakly-coupled System]{On $p$-Fractional weakly-coupled System with Critical Nonlinearities}
\author[N. Biswas and S. Chakraborty]{Nirjan Biswas$^1$\,\orcidlink{0000-0002-3528-8388} \and Souptik Chakraborty$^2$\,\orcidlink{0009-0004-1867-0560}}
\address{\rm  $^1$Department of Mathematics, Indian Institute of Science Education and Research Pune \\
Dr. Homi Bhabha Road, Pune 411008, India}
\address{\rm $^2$Tata Institute of Fundamental Research, Centre For Applicable Mathematics \\ Post Bag No 6503, GKVK Post Office, Sharada Nagar, Chikkabommasandra, Bengaluru 560065, India}
\email[N. Biswas]{nirjaniitm@gmail.com, nirjan.biswas@acads.iiserpune.ac.in}
\email[S. Chakraborty]{souptik25@tifrbng.res.in, soupchak9492@gmail.com}
\thanks{$^2$Corresponding author}

\subjclass[2020]{Primary: 35R11,  35A15, 35B33, 35J60}
\keywords{global compactness results, critical exponent problem, Palais-Smale decomposition, existence of a positive solution.}
\begin{document}
\begin{abstract}
This paper deals with the following nonlocal system of equations:
\begin{equation}
	\tag{$\mathcal S$}\label{MAT1}
	\left\{\begin{aligned}
		&(-\Delta_p)^s u = \frac{\alpha}{p_s^*}|u|^{\alpha-2}u|v|^{\beta}+f(x)\;\;\text{in}\;\mathbb{R}^{d},\\
		&(-\Delta_p)^s v = \frac{\beta}{p_s^*}|v|^{\beta-2}v|u|^{\alpha}+g(x)\;\;\text{in}\;\mathbb{R}^{d},\\
		& u, \, v >0\,  \mbox{ in }\,\mathbb{R}^{d},
	\end{aligned}
	\right.
\end{equation}
where $0<s<1<p< \infty$, $d>sp$, $\alpha,\beta>1$, $\alpha+\beta=\frac{dp}{d-sp}$, and $f,g$ are nontrivial nonnegative functionals in the dual space of $\mathcal{D}^{s,p}(\mathbb{R}^{d})$, i.e., $f,g \neq 0$, and $\prescript{}{(\mathcal{D}^{s,p})'}{\langle}f,u{\rangle}_{\mathcal{D}^{s,p}}\geq 0, \,\prescript{}{(\mathcal{D}^{s,p})'}{\langle}g,u{\rangle}_{\mathcal{D}^{s,p}}\geq 0$, whenever $u$ is a nonnegative function  in $\mathcal{D}^{s,p}(\mathbb{R}^{d})$. The primary objective of this paper is to present a global compactness result that offers a complete characterization of the Palais-Smale sequences of the energy functional associated with \eqref{MAT1}. Using this characterization, within a certain range of $s$, we establish the existence of a solution with negative energy for \eqref{MAT1} when $\ker(f)=\ker(g)$, and $\|f\|_{(\mathcal{D}^{s,p})'}, \|g\|_{(\mathcal{D}^{s,p})'}$ are small.
\end{abstract}
\maketitle
\section{Introduction}
For $0<s<1<p<\infty$ and $d>sp$, we consider the following system of equations:
\begin{equation}\tag{$\mathcal S$}\label{system_eqn}
\left\{\begin{aligned}
&(-\Delta_p)^s u = \frac{\al}{p_{s}^{*}} |u|^{\alpha-2}u|v|^{\beta} +f(x)\;\;\text{in}\;\rd,\\
&(-\Delta_p)^s v = \frac{\beta}{p^{*}_{s}} |v|^{\beta-2}v|u|^{\alpha} +g(x)\;\;\text{in}\;\rd,\\
& u, \, v >0\,  \mbox{ in }\, \rd,
\end{aligned}
\right.
\end{equation}
where $p^*_s = \frac{dp}{d-sp}$ is the critical Sobolev exponent, $\al, \be > 1, \al + \be = p^*_s$, and $f,g$ are non-negative functionals. The fractional $p$-Laplace operator $(-\Delta _p)^s$ is defined as
\begin{equation*}
   (-\Delta_{p})^{s}u(x) = 2 \lim_{\var \rightarrow 0^{+}} \int_{\mathbb{R}^{d} \backslash B(x, \var)} \frac{|u(x) - u(y)|^{p-2}(u(x)-u(y))}{|x-y|^{d+sp}}\, \dy, \; \text{for}~x \in \mathbb{R}^{d},
\end{equation*} 
where $B(x, \var)$ denotes the ball of radius $\var$ with center at $x \in \rd$. The fractional homogeneous Sobolev space $\mathcal{D}^{s,p}(\rd)$ (for brevity, we denote this space as $\wps$) is defined as the completion of $\mathcal{C}^{\infty}_{c}(\rd)$ under the Gagliardo seminorm
$$ \norm{u}_{\wps}\coloneqq \left( \iint_{\R^{2d}}\frac{|u(x)-u(y)|^p}{|x-y|^{d+sp}}\dxy \right)^{\frac{1}{p}},$$ 
which actually becomes a norm on $\wps$ as a consequence of the fractional Sobolev inequality (see \eqref{HS1}). The space $\wps$ has the following characterization (see \cite[Theorem 3.1]{Brasco2019characterisation}):
\begin{equation*}
    \wps\coloneqq\left\{u\in L^{p^*_s}(\rd):\norm{u}_{\wps}<\infty\right\},
\end{equation*}
where $\norm{\cdot}_{\wps}$ is an equivalent norm in $\wps$ and it is a reflexive Banach space.
For details on $\wps$ and its associated embedding results, we refer to \cite{Brasco2019characterisation, DiPaVa} and the references therein. We consider $\mathcal{W} = \wps \times \wps$ as the solution space for \eqref{system_eqn}, endowed with the following norm:
\begin{align*}
    \norm{(u,v)}_{\W} \coloneqq \left( \norm{u}^p_{\wps} + \norm{v}^p_{\wps} \right)^{\frac{1}{p}}, 
\end{align*}
and $f, g \in (\wps)'$, which is the dual space of $\wps$. A pair $(u,v) \in \mathcal{W}$ is a weak solution of \eqref{system_eqn}, if it satisfies the following identity for every $(\phi, \psi) \in \mathcal{W}$:
\begin{align}\label{weak1}
    &\iint\limits_{\rd\times\rd}\frac{|u(x)-u(y)|^{p-2}(u(x)-u(y))(\phi(x)-\phi(y))}{|x-y|^{d+sp}} \dxy \no \\
    &+\iint\limits_{\rd\times\rd}\frac{|v(x)-v(y)|^{p-2}(v(x)-v(y))(\psi(x)-\psi(y))}{|x-y|^{d+sp}} \dxy \no \\ 
    &=\frac{\al}{p^*_s} \int_{\rd} \abs{u}^{\al-2}u\abs{v}^{\beta} \phi \dx + \frac{\be}{p^*_s} \int_{\rd} \abs{v}^{\be-2}v\abs{u}^{\al} \psi \dx +\prescript{}{(\wps)'}{\langle}f,\phi{\rangle}_{\wps}+\prescript{}{(\wps)'}{\langle}g,\psi{\rangle}_{\wps}.
\end{align}

Elliptic systems play a crucial role in various fields, including biological applications (e.g., population dynamics) and physical systems (e.g., nuclear reactor models). These systems have garnered significant attention, as detailed in articles such as \cite{AdeMS, CFMT, Mitidieri, Reichel} and the references they cite. Additionally, problems involving the fractional Laplace and fractional $p$-Laplace operators arise in diverse areas, including phase transitions, flame propagation, chemical reactions in liquids, population dynamics, and finance. For more information, see \cite{Cafarelli, DiPaVa}.

We note that system \eqref{system_eqn} is variational, and the underlying functional is
\begin{align}\label{energy-functional}
 I_{f,g}(u,v)\coloneqq \frac{1}{p}\|(u,v)\|^p_{\mathcal{W}}
-\frac{1}{p^*_s}\int_{\rd} |u|^{\alpha}|v|^{\beta} \dx & -\prescript{}{(\wps)'}{\langle}f,u{\rangle}_{\wps} \no \\
& -\prescript{}{(\wps)'}{\langle}g,v{\rangle}_{\wps}, \; \forall \, (u,v) \in \W.
\end{align}
Observe that $I_{f,g}$ is well defined in $\W$ and $I_{f,g} \in \C^1\big(\mathcal{W})$. Moreover, if $(u,v)$ is a weak solution of~\eqref{system_eqn}, then $(u,v)$ is a critical point of $I_{f,g}$ and vice versa. 

\begin{definition}
    A sequence $\{ (u_n, v_n ) \} \subset \mathcal{W}$ is said to be a Palais-Smale (PS) sequence for $I_{f,g}$ at level $\eta$, if $I_{f, g}(u_n,v_n) \ra \eta$ in $\R$ and $I'_{f, g}(u_n,v_n) \ra 0$ in $(\mathcal{W})'$ as $n \ra \infty$. The function $I_{f,g}$ is said to satisfy (PS) condition at level $\eta$, if every (PS) sequence at level $\eta$ has a convergent subsequence.
\end{definition}

In this paper, we study the global compactness results for \eqref{system_eqn}. It is worth mentioning that every (PS) sequence of $I_{f,g}$ may not converge strongly due to the noncompactness of the embedding $\W \hookrightarrow L^{p^*_s}(\rd) \times L^{p^*_s}(\rd)$ and the presence of $\int_{\rd} \abs{u}^{\alpha} \abs{v}^{\beta}$ with $\al+\be=p^*_s$ in $I_{f,g}$. Moreover, the weak limit of the (PS) sequence can be zero even if $\eta>0$. In this paper, we classify (PS) sequence for the functional $I_{f,g}$. The classification of (PS) sequence was first studied by Struwe \cite{Struwe} for the following functional: 
\begin{align*}
    I_{\la}(u) = \frac{1}{2} \int_{\Omega} \abs{\Gr u}^2 - \frac{\la}{2} \int_{\Omega} u^2 - \frac{1}{2^*} \int_{\Omega} \abs{u}^{2^*},\; u \in \mathcal{D}_0^{1,2}(\Omega), 
\end{align*}
where $\la \in \R$, $\Omega$ is a smooth bounded domain in $\rd$ with $d>2$, $2^*=\frac{2d}{d-2}$ is the critical exponent, and $\mathcal{D}_0^{1,2}(\Omega)= \{ u \in L^{2^*}(\rd) : \abs{\Gr u} \in L^2(\rd), u=0 \text{ in } \rd \setminus \Omega \}$. Observe that every critical point of $I_{\la}$ weakly solves the Br\'ezis-Nirenberg problem 
\begin{equation}\label{brezis-nirernberg}
   -\Delta u - \la u= \abs{u}^{2^*-2}u \text{ in } \Omega; u=0 \text{ on } \pa \Omega. 
\end{equation}
Due to the presence of the term $\int_{\Omega} \abs{u}^{2^*}$, the functional $I_{\la}$ does not satisfy the (PS) condition. In \cite{Struwe}, Struwe provided a complete description of the (PS) sequence for $I_{\la}$. More precisely, the author proved that if $\{ u_n \}$ is a (PS) sequence of $I_{\la}$ at level $c$, then there exist an integer $k \ge 0$, sequences $\{ x_n^i \}_n \subset \rd, \{ r_n^i \}_n \subset \R^+$, a set of functions $u \in \mathcal{D}_0^{1,2}(\Omega), \tilde{u}_i \in \mathcal{D}^{1,2}(\rd)$ for $1 \le i \le k$ (where $\mathcal{D}^{1,2}(\rd) = \{u \in L^{2^*}(\rd) : \abs{\Gr u} \in L^2(\rd)\}$), such that $u$ weakly solves \eqref{brezis-nirernberg}, $\tilde{u}_i$ weakly solves the purely critical problem on $\rd$, i.e., $-\Delta \tilde{u}_i = \abs{\tilde{u}_i}^{2^*-2}\tilde{u}_i$ in $\rd$ such that the following hold: 
\begin{align*}
    u_n = u + \sum_{i=1}^k \tilde{u}_i^{r_n^i, x_n^i} + o_n(1), \text{ in } \mathcal{D}^{1,2}(\rd),
\end{align*}
where $o_n(1) \ra 0$ as $n \ra \infty$, and
\begin{align*}
    \tilde{u}_i^{r, y}(x) \coloneqq  r^{-\frac{d
    -2}{2}} \tilde{u}_i \left( \frac{x-y}{r} \right), \text{ for } x,y \in \rd,\,r>0.
\end{align*}
Moreover, the energy level $c$ is distributed in the following manner:
\begin{align*}
    c = I_{\la}(u) + \sum_{i=1}^k I_{\infty}(\tilde{u}_i) + o_n(1), \text{ where } I_{\infty}(u) = \frac{1}{2} \int_{\rd} \abs{\Gr u}^2 - \frac{1}{2^*} \int_{\rd} \abs{u}^{2^*}, \; u \in \mathcal{D}^{1,2}(\rd).   
\end{align*}
This result is valuable to investigate the existence of ground states in nonlinear Sch\"{o}dinger equations, Yamabe-type equations, and various types of minimization problems. In \cite[Theorem 1.2]{Willem-2010}, for $p \in (1,\infty)$, Mercuri-Willem studied a similar (PS) decomposition of signed (PS) sequences for the following functional:
\begin{align*}
    I_p(u) \coloneqq \frac{1}{p} \int_{\Omega} \abs{\Gr u}^p + \frac{1}{p} \int_{\Omega} g \abs{u}^p - \frac{\mu}{p^*} \int_{\Omega} \abs{u}^{p^*}, \; u \in \mathcal{D}_0^{1,p}(\Omega),
\end{align*}
where $\mu>0$, $\Omega$ is a smooth bounded domain in $\rd$ with $d>p$, $g \in L^{\frac{d}{p}}(\Omega)$, $p^*=\frac{dp}{d-p}$ is the critical exponent, and $\mathcal{D}_0^{1,p}(\Omega)= \{ u \in L^{p^*}(\rd) : \abs{\Gr u} \in L^p(\rd), u=0 \text{ in } \rd \setminus \Omega \}$. In \cite[Theorem 1.1]{Brasco-2016}, Brasco et. al. extended this result in the nonlocal framework by examining the global compactness property of Palais-Smale sequences associated with the following nonlocal functional:
\begin{align*}
    I_{p,s}(u) \coloneqq \frac{1}{p} \norm{u}_{\wps}^p + \frac{1}{p} \int_{\Omega} g \abs{u}^p - \frac{\mu}{p^*_s} \int_{\Omega} \abs{u}^{p^*_s}, \; u \in \mathcal{D}_0^{s,p}(\Omega),
\end{align*}
where $\Omega$ is a smooth bounded domain in $\rd$ with $d>sp$, $g \in L^{\frac{d}{sp}}(\Omega)$, and $\mathcal{D}_0^{s,p}(\Omega) = \{u \in \wps : u=0 \text{ in } \rd \setminus \Omega \}$. Further, in \cite[Theorem 1.3]{Brasco-2016}, the authors studied the global compactness result for radially symmetric functions in a ball $B \subset \rd$. For $p=2$ and $f \in (\mathcal{D}^{s,2})'$, Bhakta-Pucci in \cite[Proposition 2.1]{Bhakta-Pucci} classified the (PS) sequences associated with the following energy functional:
\begin{equation}
   I_{g,f}(u) \coloneqq \frac{1}{2} \norm{u}_{\mathcal{D}^{s,2}}^2 - \frac{1}{2^*_s} \int_{\rd} g \abs{u}^{2^*_s} -  \prescript{}{(\mathcal{D}^{s,2})'}{\langle}f,u{\rangle}_{\mathcal{D}^{s,2}},\;  u \in \mathcal{D}^{s,2},
\end{equation}
where $0<g \in L^{\infty}(\rd)$, $g(x) \ra 1$ as $\abs{x} \ra \infty$. More precisely, they established that if $\{ u_n \}$ is a (PS) sequence of $I_{g,f}$ at level $c$, then there exist an integer $k \ge 0$, sequences $\{ x_n^i \}_n \subset \rd, \{ r_n^i \}_n \subset \R^+$, a set of functions $u, \tilde{u}_i \in \mathcal{D}^{s,2}$ for $1 \le i \le k$, such that $r_n^i \ra 0, \text{ and either } x_n^i \ra x^i \in \rd \text{ or } \abs{x_n^i} \ra \infty, \text{ for } 1 \le i \le k$,
$u$ weakly solves $(-\Delta)^s u = g(x) \abs{u}^{2^*_s -2}u + f$ in $\rd$, and $\tilde{u}_i \not\equiv 0$ weakly solves the corresponding homogeneous equation $(-\Delta)^s u = g(x^i) \abs{u}^{2^*_s -2}u$ in $\rd$ such that following hold:
\begin{align*}
    u_n = u + \sum_{i=1}^k g(x^i)^{-\frac{d-2s}{4s}} \tilde{u}_i^{r_n^i, x_n^i} + o_n(1), \text{ in } \mathcal{D}^{s,2},
\end{align*}
where $o_n(1) \ra 0$ as $n \ra \infty$, and 
\begin{align*}
    \tilde{u}_i^{r, y}(x) \coloneqq  r^{-\frac{d
    -2s}{2}} \tilde{u}_i \left( \frac{x-y}{r} \right), \text{ for } x,y \in \rd,\,r>0.
\end{align*}
In this case, the energy level $c$ is distributed in the following manner:
\begin{align*}
    c = I_{g,f}(u) + \sum_{i=1}^k g(x^i)^{-\frac{d-2s}{2s}} I_{1,0}(\tilde{u}_i) + o_n(1).   
\end{align*}
Observe that, by the uniqueness of the positive solution of 
\begin{align*}
    (-\Delta)^s u = u^{2^*_s-1} \text{ in } \rd,\, u \in \D^{s,2}.
\end{align*}
for each $i$, $\tilde{u}_i$ is a nonlocal Aubin-Talenti bubble (up to translation and dilation). In \cite[Proposition 2.1]{Bhakta-Pucci} the authors established the following bubble interaction 
\begin{align}\label{bub-1}
\left| \log \left( \frac{r_n^i}{r_n^j} \right) \right| + \left| \frac{x_n^i - x_n^j}{r_n^i}  \right| \rightarrow \infty, \text{ for } 1 \le i \le k.
\end{align}
The above interaction in \eqref{bub-1} is motivated by the works of Palatucci-Pisante (see \cite{GiAd2014, Palatucci-Pisante-NA}). In \cite{GiAd2014}, Palatucci-Pisante measured the non-compactness of the embedding $\mathcal{D}^{s,2} \hookrightarrow L^{2_s^*}(\rd)$ under the conformal group action of translation and dilation. Using the profile decomposition for general bounded sequence on $\mathcal{D}^{s,2}$ and asymptotic orthogonal relation between the weakly interacting bubbles (see \cite{GiAd2014}), they proved (see \cite[Theorem~1.2]{Palatucci-Pisante-NA}) the profile decomposition for the (PS) sequence of the energy functional associated with fractional Br\'ezis Nirenberg problem (see \cite{Ser-Val}).

The classification of (PS) sequences associated with scalar equations (local/nonlocal) has been extensively studied in the literature, with significant contributions from several authors, including \cite{Struwe, Gerard, Palatucci-Pisante-NA, Bhakta-Pucci, Brasco-2016}, to name just a few. Smets \cite{Smets-TAMS} observed that, in the presence of Hardy potential for the local case, noncompactness arises due to concentration occurring through two distinct profiles (see also \cite{Bhak-San} for the local case involving Hardy-Sobolev-Maz'ya type equations, and \cite{BCP} for the nonlocal case with a Hardy-Sobolev term). Additionally, \cite{Tin-Fi} presents a more abstract approach for profile decomposition in general Hilbert spaces. However, to the best of our knowledge, the (PS) decomposition in the context of systems of equations has not been extensively explored. In this direction, the work of Peng-Peng-Wang \cite{PPW} investigates the (PS) decomposition for systems of equations on bounded domains in the local case when $p=2$. Subsequently, Bhakta et al. \cite{MoSoMiPa} examined the global compactness results of \eqref{system_eqn} for $p=2$.

The following theorem classifies the (PS) sequence for the functional $I_{f,g}$ defined in \eqref{energy-functional}. 

\begin{theorem}[(PS) decomposition for nonhomogeneous system]\label{PS-decomposition}
   Let $s\in (0,1)$ and $p \in (1, \infty)$. Let $\al, \be > 1,\, \al + \be = p^*_s$.  Assume that  $f, g$ are nontrivial nonnegative functionals in the dual space of $\wps$. Let $\{ (u_n, v_n ) \}$ be a (PS) sequence for $I_{f,g}$ at level $\eta$. Then there exists a subsequence (still denoted by $\{ (u_n, v_n ) \}$) for which the following hold: 
    \noi there exist an integer $k \ge 0$, sequences $\{ x_n^i \}_n \subset \rd, \{ r_n^i \}_n \subset \R^+$, pair of functions $(u, v), (\tilde{u}_i, \tilde{v}_i) \in \W$ for $1 \le i \le k$ such that $(u, v)$ satisfies \eqref{system_eqn} without sign assumptions and $(\tilde{u}_i, \tilde{v}_i)$ satisfies the following homogeneous system \eqref{system_eqn} with $f=g=0$:
    \begin{equation}\label{homogenous-system}
\left\{\begin{aligned}
&(-\Delta_p)^s \tilde{u}_i = \frac{\al}{p_{s}^{*}} |\tilde{u}_i|^{\alpha-2} \tilde{u}_i |\tilde{v}_i|^{\beta}\;\;\text{in}\;\rd,\\
&(-\Delta_p)^s \tilde{v}_i = \frac{\beta}{p^{*}_{s}} |\tilde{v}_i|^{\beta-2} \tilde{v}_i |\tilde{u}_i|^{\alpha}\;\;\text{in}\;\rd,
\end{aligned}
\right.
\end{equation}
such that 
\begin{align*}
    &(u_n ,v_n) = (u, v) + \sum_{i=1}^k (C_{x_n^i, r_n^i}\tilde{u}_i, C_{x_n^i, r_n^i}\tilde{v}_i) + o_n(1), \text{ and } \eta = I_{f,g}(u,v) + \sum_{i=1}^k I_{0,0}(\tilde{u}_i , \tilde{v}_i) + o_n(1), \\
    &\text{either } x_n^i \to x^i, \text{and } r_n^i \ra 0, \text{or } \infty, \text{ or } \abs{x_n^i} \ra \infty, \text{ for } 1 \le i \le k, \\
    & \left| \log \left( \frac{r_n^i}{r_n^j} \right) \right| + \left| \frac{x_n^i - x_n^j}{r_n^i}  \right| \rightarrow \infty, \text{ for } i \neq j, 1 \le i, j \le k,
\end{align*}
where $o_n(1) \ra 0$ in $\mathcal{W}$ and $C_{x_n^i, r_n^i}\tilde u_i(x)\coloneqq (r_n^i)^{-\frac{d-sp}{p}}\tilde u_i(\frac{x-x_n^i}{r_n^i})$ and in the case $k=0$, the above expression holds without  $(\tilde{u}_i, \tilde{v}_i),\, x_n^i$ and $r_n^i$. 
\end{theorem}

Using the above (PS) decomposition result, we prove the existence of a positive solution of \eqref{system_eqn}, in the spirit of \cite{MoSoMiPa} (also see \cite{Jeanjean, Bhakta-Pucci}). There is a limited number of work on weakly coupled systems of equations in the nonlocal case. Notable contributions in this area include \cite{PPW, Chen-Squassina, Faria}, which studies Dirichlet systems of equations on bounded domains. For nonlocal nonhomogeneous fractional Laplacian and fractional $p$-Laplacian systems on the entire $\rd$, we refer to \cite{MoSoMiPa, BCP-Sys, Fiscella-Pucci-Zhang, Mukherjee-Sreenadh} and the references therein. To the best of our knowledge, no prior work has addressed the existence of positive solutions to \eqref{system_eqn} involving the fractional $p$-Laplacian and critical exponents on $\rd$. 


\begin{theorem}\label{th:mul}
Let $p \in (1, \infty)$ and $s \in(0, \frac{p-1}{p})$. Let $\al, \be > 1,\, \al + \be = p^*_s$.  Assume that  $f, g$ are  nontrivial nonnegative functionals in the dual space of $\wps$ with $ker(f)=ker(g)$ and
$$ \max\{\|f\|_{(\wps)'}, \|g\|_{(\wps)'}\} \le \frac{1}{2} \left( \frac{1}{p} - \frac{p-1}{p^*_s(p^*_s-1)} \right) \left(  \frac{p-1}{p^*_s-1} S_{\al, \be}^{\frac{p^*_s}{p}} \right)^{\frac{p-1}{p^*_s-p}},$$
where $S_{\al, \be}$ is the best constant of \eqref{HS2},
then \eqref{system_eqn} admits a positive solution with negative energy. Furthermore, if $f\equiv g$, then the solution $( u, v)$ of \eqref{system_eqn} has the property that
$ u\not\equiv v$, whenever $\alpha\neq\beta$. Finally, if $\alpha=\beta$ but $f\not\equiv g$, then $u\not\equiv  v$.
\end{theorem}

To prove Theorem~\ref{th:mul}, we show that an auxiliary functional $J_{f,g}$ (defined in \eqref{20-24-5'} and associated with $I_{f,g}$) attains its infimum on the exterior of a ball in $\W$. The existence of this solution is based on the following statement
\begin{align}\label{converse-1}
    \left|\log \left( \frac{r_n^i}{r_n^j} \right)\right| + \left| \frac{x_n^i - x_n^j}{r_n^i} \right| \rightarrow \infty \Longrightarrow \mathcal{A}(C_{x_n^i,r_n^i}u,C_{x_n^j,r_n^j}v)\to 0,
\end{align}
where $i \neq j, 1 \le i, j \le k$ and the map $\mathcal{A}$ is defined as
\begin{align*}
    \mathcal{A}(u,v) \coloneqq \iint\limits_{\rd \times \rd} \frac{\abs{u(x) - u(y)}^{p-2} (u(x) - u(y)) (v(x) - v(y))}{\abs{x-y}^{d+sp}} \dxy, \; \forall \, (u,v) \in \mathcal{W}.
\end{align*}

In simple terms, \eqref{converse-1} says two arbitrary $\wps$-functions are considered asymptotically orthogonal (in this case, this means, for every $u,\,v\in\wps$, $\mathcal{A}(C_{x_n^i,r_n^i}u, C_{x_n^j,r_n^j}v)\to 0$ as $n\to \infty$) when sequences of group elements acting on them go to infinity in different directions. As a result, we will observe that asymptotic orthogonality is achieved either when the dilation parameters are not comparable or when they are comparable but negligible relative to the distance between the translation parameters as $n\to\infty$. In \cite[Lemma~3]{GiAd2014}, for $p= 2$, the authors established \eqref{converse-1} using stronger density results for $\mathcal{D}^{s,2}$ which mainly depends on the Fourier representation of the norm in $\mathcal{D}^{s,2}$ and the operator $(-\Delta)^s$. However, adapting those techniques to the case $p \neq 2$ is not straightforward. For $p \in (1, \infty)$ and $s \in (0, \frac{p-1}{p})$, in Proposition \ref{weak-bub-converse-II} we show that \eqref{converse-1} holds, and our approach relies on this restrictive range of $s$. Now, a natural question to ask: 

\noi \textbf{Question:}
Whether \eqref{converse-1} holds for every $s \in (0,1)$ and $p \neq 2$?

The structure of the paper is as follows: In Section 2, we present the convergence of certain integrals, the behaviour of bounded sequences in $\wps$ under the action of a group of isometries, refined embeddings of $\wps$ into the Morrey spaces, and the global compactness results for \eqref{system_eqn}. In this section, we prove Theorem \ref{PS-decomposition}. In the last section, we prove Theorem \ref{th:mul}. 

\medskip
\noi \textbf{Notation:} Throughout this paper, we use the following notation and definition:

\noi (a) $B(x,r)$ denotes a ball of radius $r$ with center at $x \in \rd$. $B(x,r)^c := \rd \setminus B(x,r)$, and $\omega_N := |B(x,1)|$.

\noi (b) $p'\coloneqq\frac{p}{p-1}$ denotes the H\"{o}lder's conjugate exponent of $p>1$.

\noi (c) For a function $f$, the positive part is defined as $f_+(x):= \max \{f(x), 0\}$ and the negative part is defined as $f_-(x):= -\min \{f(x), 0\}$.

\noi (d) $\mathcal{S}(\rd)$ denotes the Schwartz space.

\noi (e) $P\lesssim_{a}Q$ represents that there exists a constant $C\equiv C(a)>0$ such that $P\leq CQ$ holds.
\section{Global compactness results}
\subsection{Technical lemmas}
This section examines the convergence of certain integrals that will be useful in the subsequent sections. First, we recall the classical fractional Sobolev inequality (see \cite[Theorem 2.2.1]{BV}):
\begin{equation}\label{HS1}
\left( \int_{\rd} \abs{u(x)}^{p^*_s} \dx \right)^{\frac{p}{p^*_s}} \lesssim_{d,p,s}  \iint\limits_{\rd \times \rd}\frac{|u(x)-u(y)|^p}{|x-y|^{d+sp}}\dxy, \; \forall \, u \in \wps. 
\end{equation}
Using \eqref{HS1}, the following inequality holds for $\al, \be \ge 1,\, \al + \be = p^*_s$:
\begin{align}\label{HS2}
    &\int_{\rd} \abs{u(x)}^{\al} \abs{v(x)}^{\be} \dx \le \left( \int_{\rd} \abs{u(x)}^{p^*_s} \dx \right)^{\frac{\al}{p^*_s}} \left( \int_{\rd} \abs{v(x)}^{p^*_s} \dx \right)^{\frac{\be}{p^*_s}} \no \\
    &\lesssim_{d,p,s} \left( \; \iint\limits_{\rd \times \rd}\frac{|u(x)-u(y)|^p}{|x-y|^{d+sp}}\dxy \right)^{\frac{\al}{p}} \left( \; \iint\limits_{\rd \times \rd}\frac{|v(x)-v(y)|^p}{|x-y|^{d+sp}}\dxy \right)^{\frac{\be}{p}}, \no \\
    &\lesssim_{d,p,s} \norm{(u,v)}_{\W}^{\frac{p^*_s}{p}}, \; \forall \, (u,v) \in \W.
\end{align}

We recall the following lemma from \cite[Theorem 2]{brezis1983relation}.

\begin{lemma}\label{convergence1}
    Let $j : \mathbb{C} \ra \mathbb{C}$ be a continuous function with $j(0) = 0$. Given $\var>0$, assume that there exists two continuous functions $\phi_{\var}$ and $\psi_{\var}$ such that the following inequality holds: 
\begin{align}\label{condition-j}
    \abs{j(a+b) - j(a)} \le \var \phi_{\var}(a) + \psi_{\var}(b), \; \forall \, a,b \in \mathbb{C}. 
\end{align}
Further, let $f_n: \rd \ra \mathbb{C}$ with $f_n = f+ g_n$ be a sequence of measurable functions such that 
\begin{enumerate}
    \item[\rm{(i)}] $g_n \ra 0$ a.e. in $\rd$,
    \item[\rm{(ii)}] $j(f) \in L^1(\rd)$,
    \item[\rm{(iii)}] $\int_{\rd} \phi_{\var}(g_n(x)) \dx \le C < \infty$, for some constant $C$ independent of $\var$ and $N$. 
    \item[\rm{(iv)}] $\int_{\rd} \psi_{\var}(f(x)) \dx < \infty$ for all $\var> 0$.
\end{enumerate}
Then 
\begin{align*}
   \lim_{n \ra \infty} \int_{\rd} \left| j(f_n) - j(f) - j(g_n) \right| \dx = 0.
\end{align*}
\end{lemma}

The above lemma yields the following convergence.

\begin{lemma}\label{convergence2}
Let $\al, \be > 1$ and $\al + \be = p^*_s$. If $u_n \rightharpoonup u$ and $v_n \rightharpoonup v$ in $\wps$ with $u_n(x)\to u(x)$ and $v_n(x)\to v(x)$ a.e. $x \in \rd$, then 
\begin{align*}
   \lim_{n \ra \infty} \int_{\rd} \left( \abs{u_n}^{\alpha} \abs{v_n}^{\beta} - \abs{u}^{\alpha} \abs{v}^{\beta} - \abs{u_n - u}^{\alpha} \abs{v_n-v}^{\beta} \right) \dx = 0.
\end{align*}
\end{lemma}

\begin{proof}
    Consider $j : \R^2 \ra \R$ such that $j(x,y) = \abs{x}^{\al} \abs{y}^{\be}$. Then, in view of the following inequality:
    \begin{align*}
        \left| \abs{x+a}^{\alpha} \abs{y+b}^{\beta} - \abs{x}^{\al} \abs{y}^{\be} \right| \le \var \left( \abs{x}^{p^*_s} + \abs{y}^{p^*_s} \right) + C_{\var} \left(\abs{a}^{p^*_s} + \abs{b}^{p^*_s} \right), \; \forall \, x,y,a,b \in \mathbb{R}, 
    \end{align*}
    which holds similarly as in \cite[Lemma 3.2]{MoSoMiPa}, $j$ satisfies \eqref{condition-j} with $\phi_{\var}(x,y) = \psi_{\var}(x,y)= \abs{x}^{p^*_s} + \abs{y}^{p^*_s}$.  Now we consider 
    \begin{align*}
        f_n \coloneqq (u_n, v_n), f \coloneqq (u,v), \text{ and } g_n \coloneqq (u_n - u, v_n -v). 
    \end{align*}
    Using the boundedness of $\{ u_n \}, \{ v_n \}$ in $\wps$, 
    \eqref{HS1} and \eqref{HS2}, we can verify that (i)-(iv) hold in Lemma \ref{convergence1}. Therefore, we get the required convergence by applying Lemma \ref{convergence1}.   
\end{proof}

Now we discuss the classical Brézis–Lieb Lemma and its consequences. 

\begin{lemma}\label{BL}
    Let $1<q< \infty$. Let $\{ f_n \} \subset L^{q}(\rd)$ be a bounded sequence such that $f_n(x) \ra f(x)$ a.e. $x \in \rd$. Then the following hold: 
    \begin{enumerate}
        \item[\rm{(i)}] $\norm{f_n}_{L^q(\rd)}^q - \norm{f_n - f}_{L^q(\rd)}^q + o_n(1) = \norm{f}_{L^q(\rd)}^q.$
        \item[\rm{(ii)}] Consider the function $J_q$ defined as $J_q(t)=\abs{t}^{q-2}t$. Then 
        \begin{align*}
           J_q(f_n) - J_q(f_n-f) = J_q(f) + o_n(1) \text{ in } L^{q'}(\rd).
        \end{align*}
    \end{enumerate}
\end{lemma}
\begin{proof}
Proof of (i) follows from \cite{Br-Li}, and proof of (ii) follows from \cite[Lemma 3.2]{Willem-2010}. 
\end{proof}
The above lemma leads to the following convergence. A similar convergence result has been proved in \cite{Brasco-2016}. 

\begin{lemma}\label{convergence-BL}
    Let $\{ (u_n, v_n) \}$ weakly converge to $(u,v)$ in $\W$ with $u_n(x) \ra u(x)$ and $v_n(x) \ra v(x)$ a.e. $x \in \rd$. Then, up to a subsequence, the following hold
   \begin{enumerate}
       \item[\rm{(i)}] $\norm{u_n}_{\wps}^p - \norm{u_n - u}_{\wps}^p = \norm{u}^p_{\wps} + o_n(1), \text{ and } \norm{v_n}_{\wps}^p - \norm{v_n - v}_{\wps}^p = \norm{v}_{\wps}^p + o_n(1).$ 
       \item[\rm{(ii)}] Consider the function $J_p$ defined as $J_p(t)=\abs{t}^{p-2}t$. Then 
       \begin{align*}
           \frac{J_p(u_n(x) - u_n(y))}{\abs{x-y}^{\frac{d+sp}{p'}}} - \frac{J_p((u_n(x)-u(x)) - (u_n(y)-u(y))}{\abs{x-y}^{\frac{d+sp}{p'}}} = \frac{J_p(u(x) - u(y))}{\abs{x-y}^{\frac{d+sp}{p'}}} + o_n(1),
       \end{align*}
       in $L^{p'}(\mathbb{R}^{2d})$.
   \end{enumerate}

\end{lemma}
\begin{proof}
(i) We choose 
    \begin{align*}
        f_n(x,y) = \frac{u_n(x) - u_n(y)}{ \abs{x-y}^{\frac{d+sp}{p}} },\, f(x,y) = \frac{u(x) - u(y)}{ \abs{x-y}^{\frac{d+sp}{p}} }, \text{ and } N=2d,
    \end{align*}
     in Lemma \ref{BL} to get $\norm{u_n}_{\wps}^p - \norm{u_n - u}_{\wps}^p = \norm{u}^p_{\wps} + o_n(1)$. Similarly, $\norm{v_n}_{\wps}^p - \norm{v_n - v}_{\wps}^p = \norm{v}_{\wps}^p + o_n(1)$.

\noi (ii) Proof follows with the same choices of $f_n,\,f$ and using (ii) of Lemma \ref{BL}.
\end{proof}

The following lemma states the convergence of some integrals.

\begin{lemma}\label{convergence-integrals}
    Let $\al, \be > 1$ and $\al + \be = p^*_s$. Let $\{ (u_n, v_n) \}$ weakly converge to $(u,v)$ in $\W$. 
    \begin{enumerate}
        \item[\rm{(i)}] Then up to a subsequence
    \begin{align*}
        \lim_{n \ra \infty} \int_{\rd} \abs{u_n}^{\al -2} u_n \abs{v_n}^{\beta} \phi = \int_{\rd} \abs{u}^{\al -2} u \abs{v}^{\beta} \phi, \text{ and }  \lim_{n \ra \infty} \int_{\rd} \abs{u_n}^{\al} \abs{v_n}^{\be -2} v_n \psi = \int_{\rd} \abs{u}^{\al} \abs{v}^{\be -2} v \psi,
    \end{align*}
    for all $(\phi, \psi) \in \W$.
    \item[\rm{(ii)}] Then up to a subsequence
    \begin{align*}
    &\lim_{n \ra \infty} \iint\limits_{\rd\times\rd}\frac{|u_n(x)-u_n(y)|^{p-2}(u_n(x)-u_n(y))(\phi(x)-\phi(y))}{|x-y|^{d+sp}} \dxy  \\
    &=\iint\limits_{\rd\times\rd}\frac{|u(x)-u(y)|^{p-2}(u(x)-u(y))(\phi(x)-\phi(y))}{|x-y|^{d+sp}} \dxy, \text{ and } \\
    &\lim_{n \ra \infty} \iint\limits_{\rd\times\rd}\frac{|v_n(x)-v_n(y)|^{p-2}(v_n(x)-v_n(y))(\psi(x)-\psi(y))}{|x-y|^{d+sp}} \dxy  \\
    &=\iint\limits_{\rd\times\rd}\frac{|v(x)-v(y)|^{p-2}(v(x)-v(y))(\psi(x)-\psi(y))}{|x-y|^{d+sp}} \dxy,
    \end{align*}
    for all $(\phi, \psi) \in \W$.
    \item[\rm{(iii)}] Then up to a subsequence
    \begin{align*}
    &\lim_{n \ra \infty} \iint\limits_{\rd\times\rd}\frac{|\phi(x)-\phi(y)|^{p-2}(\phi(x)-\phi(y))(u_n(x)-u_n(y))}{|x-y|^{d+sp}} \dxy  \\
    &=\iint\limits_{\rd\times\rd}\frac{|\phi(x)-\phi(y)|^{p-2}(\phi(x)-\phi(y))(u(x)-u(y))}{|x-y|^{d+sp}} \dxy, \text{ and } \\
    &\lim_{n \ra \infty} \iint\limits_{\rd\times\rd}\frac{|\psi(x)-\psi(y)|^{p-2}(\psi(x)-\psi(y))(v_n(x)-v_n(y))}{|x-y|^{d+sp}} \dxy  \\
    &=\iint\limits_{\rd\times\rd}\frac{|\psi(x)-\psi(y)|^{p-2}(\psi(x)-\psi(y))(v(x)-v(y))}{|x-y|^{d+sp}} \dxy,
    \end{align*}
    for all $(\phi, \psi) \in \W$.
    \end{enumerate}
\end{lemma}
\begin{proof}
(i) Consider the following quantity
\begin{align*}
    M_1 = \max \left\{ \sup_{n\in\N}\norm{u_n}_{L^{p^*_s}(\rd)}^{\al -1}, \sup_{n\in\N}\norm{v_n}_{L^{p^*_s}(\rd)}^{\be}, \norm{u}_{L^{p^*_s}(\rd)}^{\al-1},  \norm{v}_{L^{p^*_s}(\rd)}^{\be} \right\}.
\end{align*}
Clearly, $M_1< \infty$ since $\{ u_n\}, \{ v_n \}$ are bounded on $\wps$. Let $\var > 0$ be given.  Take $\phi \in \wps$, and we choose $R(\var)> > 1$ such that 
\begin{align*}
   \left( \int_{\rd \setminus B(0, R(\var))} \abs{ \phi }^{p^*_s}\,\dx \right)^{\frac{1}{p^*_s}} < \frac{\var}{2M_1^2}.
\end{align*}
Now we split
\begin{align*}
    &\left| \int_{\rd} \left( \abs{u_n}^{\al -2} u_n \abs{v_n}^{\beta} - \abs{u}^{\al -2} u \abs{v}^{\beta} \right) \phi\,\dx \right| \\
    &= \left| \left( \int_{B(0, R(\var))} + \int_{\rd \setminus B(0, R(\var))} \right) \left( \abs{u_n}^{\al -2} u_n \abs{v_n}^{\beta} - \abs{u}^{\al -2} u \abs{v}^{\beta} \right) \phi\,\dx \right|. 
\end{align*}
Using H\"{o}lder's inequality with the conjugate triplet $(\frac{p^*_s}{\al-1}, \frac{p^*_s}{\be}, p^*_s)$, we calculate 
    \begin{align}\label{converge-1}
        &\left| \int_{\rd \setminus B(0, R(\var))} \left( \abs{u_n}^{\al -2} u_n \abs{v_n}^{\beta} - \abs{u}^{\al -2} u \abs{v}^{\beta} \right) \phi\,\dx \right| \le \int_{\rd \setminus B(0, R(\var))} \left( \abs{u_n}^{\al-1} \abs{v_n}^{\beta}  + \abs{u}^{\al-1} \abs{v}^{\beta} \right) \abs{\phi}\,\dx \no \\
        & \le \left( \norm{u_n}_{L^{p^*_s}(\rd)}^{\al-1} \norm{v_n}_{L^{p^*_s}(\rd)}^{\be} +  \norm{u}_{L^{p^*_s}(\rd)}^{\al-1} \norm{v}_{L^{p^*_s}(\rd)}^{\be} \right) \norm{\phi}_{L^{p^*_s}(\rd \setminus B(0, R(\var)))} < \var. 
    \end{align}
Moreover, using the compact embeddings of $\wps \hookrightarrow L_{loc}^p(\rd)$, we have $u_n \ra u$ and $v_n \ra v$ in $L^p(K)$ for any compact set $K \subset \rd$. As a consequence, up to a subsequence, we get $u_n(x) \ra u(x)$ and $v_n(x) \ra v(x)$ for a.e. $x \in \rd$. 
Hence, using the Vitali convergence theorem, we conclude that 
\begin{align}\label{converge-1.1}
 \lim_{n \ra \infty} \int_{B(0, R(\var))}  \left( \abs{u_n}^{\al -2} u_n \abs{v_n}^{\beta} - \abs{u}^{\al -2} u \abs{v}^{\beta} \right) \phi\,\dx  = 0.
\end{align}
Therefore, from \eqref{converge-1} and \eqref{converge-1.1}, we get
\begin{align*}
    \lim_{n \ra \infty} \int_{\rd}  \left( \abs{u_n}^{\al -2} u_n \abs{v_n}^{\beta} - \abs{u}^{\al -2} u \abs{v}^{\beta} \right) \phi\,\dx  = 0
\end{align*}
The second convergence also holds using a similar set of arguments. 

\noi (ii) Set $M_2 = \max \left\{ \sup_{n\in\N} \norm{u_n}_{\wps}^{p-1}, \norm{u}_{\wps}^{p-1} \right\}.$ Let $\var > 0$ be given. We split 
\begin{align*}
    \left| \iint_{\mathbb{R}^{2d}} f_n(x,y) \frac{\dxy}{\abs{x-y}^{d+sp}} \right| = \left| \left( \iint_{ B(0, R(\var))} + \iint_{\mathbb{R}^{2d} \setminus B(0, R(\var))} \right) f_n(x,y) \frac{\dxy}{\abs{x-y}^{d+sp}} \right|
\end{align*}
where $B(0, R(\var)) \subset \R^{2d}$. Take $\phi \in \wps$ and $f(x,y) = \abs{\phi(x) - \phi(y)}^p \abs{x-y}^{-(d+sp)}$ for $x,y \in \rd$. Then $\norm{f}_{L^1(\mathbb{R}^{2d})} = \norm{\phi}_{\wps}^p< \infty$. So we can choose $R(\var)> > 1$ such that 
\begin{align}\label{es-1}
   \left( \iint_{\mathbb{R}^{2d} \setminus B(0, R(\var))} \abs{\phi(x) - \phi(y)}^p \frac{\dxy}{\abs{x-y}^{d+sp}} \right)^{\frac{1}{p}} = \norm{f}_{L^1(\mathbb{R}^{2d} \setminus B(0, R(\var)))}^{\frac{1}{p}} < \frac{\var}{2M_2^2}.
\end{align}
Set $$f_n(x,y) = \left( \abs{u_n(x) - u_n(y)}^{p-2} (u_n(x) - u_n(y)) - \abs{u(x) - u(y)}^{p-2} (u(x) - u(y))\right) (\phi(x) - \phi(y))$$ for $x, y \in \rd$. Then H\"{o}lder's inequality with the conjugate pair $(p, p')$ yields
\begin{align*}
    &\iint_{\mathbb{R}^{2d} \setminus B(0, R(\var))} \frac{\abs{f_n(x,y)}}{\abs{x-y}^{d+sp}} \dxy \\
    & \le  \iint_{\mathbb{R}^{2d} \setminus B(0, R(\var))} \frac{ \abs{u_n(x) - u_n(y)}^{p-1} + \abs{u(x) - u(y)}^{p-1} }{\abs{x-y}^{d+sp}} \abs{\phi(x) - \phi(y)} \dxy \\
    & \le \left( \norm{u_n}_{\wps}^{p-1} + \norm{u}_{\wps}^{p-1} \right) \left( \iint_{\mathbb{R}^{2d} \setminus B(0, R(\var))} \abs{\phi(x) - \phi(y)}^p \frac{\dxy}{\abs{x-y}^{d+sp}} \right)^{\frac{1}{p}} < \var,
\end{align*}
where the last inequality holds using \eqref{es-1}. 
The rest of the proof follows using the same arguments as in (i). The second convergence holds similarly.

\noi (iii) Proof follows with analogous arguments as in (ii) with appropriate adjustments.  
\end{proof}

For brevity, from now onwards we denote $\mathcal{A}(u,v)$ as
\begin{align}\label{action-1}
    \mathcal{A}(u,v) \coloneqq \iint\limits_{\rd \times \rd} \frac{\abs{u(x) - u(y)}^{p-2} (u(x) - u(y)) (v(x) - v(y))}{\abs{x-y}^{d+sp}} \dxy, \; \forall \, (u,v) \in \mathcal{W}.
\end{align}

\subsection{Behavior of bounded sequences in \texorpdfstring{$\wps$}{} under the action of group of isometries}

Let $\mathcal{T},\,\mathcal{D}\subset \mathcal{U}(\wps)$ be the classes of isometric operators induced by translations and dilations, respectively on $\rd$, where
\begin{align*}
\mathcal{U}(\wps) &\coloneqq \left\{\phi:\wps\to\wps : \|\phi(u)\|_{\wps}=\|u\|_{\wps}\right\},\\
  \mathcal{T}&\coloneqq \left\{T_y\in \mathcal{U}(\wps) : y\in\rd \text{ and }T_yu(x)=u(x-y)\right\},\\
  \mathcal{D} &\coloneqq \left\{D_{\lambda}\in \mathcal{U}(\wps) : \lambda\in (0,\infty) \text{ and }D_{\lambda}u(x)=\lambda^{-\frac{d-sp}{p}}u ( \frac{x}{\lambda} ) \right\}.
\end{align*}
Now we consider the family of operators $\mathcal{C}$ given by the composition of translations and dilations on $\rd$ (actually the semidirect product $\mathcal{T}\rtimes\mathcal{D}$); i.e., we define $\mathcal{C}\subset \mathcal{U}(\wps)$ as
  \bee\no
  \mathcal{C}\coloneqq \left\{C_{y,\lambda} \in \mathcal{U}(\wps): C_{y,\lambda}u(x)\coloneqq \lambda^{-\frac{d-sp}{p}}u \left(\frac{x-y}{\lambda}\right), \; \forall \, u\in\wps;\,y\in\rd;\,\lambda>0\right\}.
  \eee 

  Consider $G=\rd\rtimes (0,\infty)$. To give a group structure on $G$, we identify $G$ with 
\begin{align*}
Aff(\rd) \coloneqq \left\{\phi_{y,\lambda}:\rd\to\rd : \phi_{y,\lambda}(x)= \frac{x-y}{\lambda}, (y,\lambda)\in G \right\}.
\end{align*}
Now,
\begin{align*}
        (\phi_{y,\lambda}\circ \phi_{a,\delta})(x) = \phi_{y,\lambda} \left(\frac{x-a}{\delta}\right) = \frac{x-(a+\delta y)}{\lambda\delta}=\phi_{a+\delta y,\lambda\delta}(x),\text{ for } x\in\rd.
\end{align*}
Define the group law on $G$ by 
\begin{align*}
(a,\delta) \circ  (y,\lambda) \coloneqq (a+\delta y,\lambda\delta).
\end{align*}
Observe that, $G$ with the above group law is nonabelian. Next, consider the following metric on $G$:
\begin{align}\label{metric-new}
     \displaystyle d\left((y,\lambda),(w,\sigma)\right) \coloneqq \left| \log\left(\frac{\lambda}{\sigma}\right) \right| + |y-w|.
\end{align}
We can verify that $G$ is complete with respect to $d$ (note that $(0,\infty)\subset \mathbb{R}$ with respect to the induced metric of $\mathbb{R}$ is not complete). The following lemma shows how the map $\mathcal{A}$ (defined in \eqref{action-1}) behaves under the above group law.
\begin{lemma}\label{Asy-Bub-inter}
Let $(a,\delta), (b, \lambda) \in \rd\rtimes (0,\infty)$. Then for every $u, v \in \wps(\rd)$ we have 
\begin{align}\label{A-behaviour}
\mathcal{A} \left( C_{a,\delta}u, C_{b, \lambda} v \right) = \mathcal{A}\left( u, C_{\frac{b-a}{\delta}, \frac{\lambda}{\delta}} v \right) = \mathcal{A}(C_{\frac{a-b}{\la},\frac{\delta}{\la}}u,v).
\end{align}
\end{lemma}
\begin{proof}
By the definition of $\mathcal{A}$,  
\begin{align*}
&\mathcal{A} \left( C_{a,\delta}u, C_{b, \lambda} v \right) \\
&= \delta^{-\frac{d-sp}{p'}} \lambda^{-\frac{d-sp}{p}} \iint\limits_{\rd \times \rd} \frac{\left| u (\frac{x-a}{\delta}) - u(\frac{y-a}{\delta}) \right|^{p-2} \left( u(\frac{x-a}{\delta}) - u(\frac{y-a}{\delta}) \right)\left( v(\frac{x-b}{\lambda}) - v(\frac{y-b}{\lambda})\right)}{\abs{x-y}^{d+sp}}\dxy\\
& =  \delta^{\frac{d-sp-2dp}{p}} \lambda^{-\frac{d-sp}{p}} \iint\limits_{\rd \times \rd} \frac{\left| u(\frac{x-a}{\delta}) - u(\frac{y-a}{\delta}) \right|^{p-2} \left( u(\frac{x-a}{\delta}) - u(\frac{y-a}{\delta}) \right)\left(v(\frac{x-b}{\lambda}) - v(\frac{y-b}{\lambda})\right)}{\abs{\frac{x-a}{\delta}-\frac{y-a}{\delta}}^{d+sp}}\dx\dy.
\end{align*}
Using the change of variable $\overline{x} = \frac{x-a}{\delta}, \overline{y} = \frac{y-a}{\delta}$, from the above identities we obtain
\begin{align*}
\mathcal{A} \left( C_{a,\delta}u, C_{b, \lambda} v \right) & = \left( \frac{\delta}{\lambda}\right)^{\frac{d-sp}{p}} \iint\limits_{\rd \times \rd} \frac{\abs{u(\overline{x}) - u(\overline{y})}^{p-2} \left( u(\overline{x}) - u(\overline{y}) \right)\left( v(\frac{\overline{x}\delta + a - b}{\lambda}) - v(\frac{\overline{y}\delta + a - b}{\lambda}) \right)}{\abs{\overline{x}-\overline{y}}^{d+sp}} \, \d \overline{x} \d \overline{y} \\
& = \mathcal{A}\left( u, C_{\frac{b-a}{\delta}, \frac{\lambda}{\delta}} v \right).
\end{align*}
Using the change of variable $\overline{x} = \frac{x-b}{\la}, \overline{y} = \frac{y-b}{\la}$, we also get $ \mathcal{A} \left( C_{a,\delta}u, C_{b, \lambda} v \right) = \mathcal{A}(C_{\frac{a-b}{\la},\frac{\delta}{\la}}u,v)$. Thus, \eqref{A-behaviour} holds. 
\end{proof}

Now, we prove that the map $G \ra \mathcal{U}(\wps)$ defined by $(y,\lambda) \mapsto C_{y,\lambda}$, is continuous with respect to the strong topology of $\mathcal{U}(\wps)$ and the metric topology on $G$ induced by $d$ defined in \eqref{metric-new}.

\begin{proposition}\label{Cylambda}
    If $(y_n,\lambda_n)\to (y,\lambda)$ in $\rd\times (0,\infty)$ with respect to the metric $d$ defined in \eqref{metric-new}, then for all $u\in \wps$, $C_{y_n,\lambda_n}u \rightarrow C_{y,\lambda}u$ in the norm topology of $\wps$.
\end{proposition}
\begin{proof}
Let $|y_n-y|+|\log(\frac{\lambda_n}{\lambda})|\rightarrow 0$ as $n\to\infty$. Then $y_n \ra y$, and moreover using $|\log(\la_n) - \log(\la)|=|\log(\frac{\lambda_n}{\lambda})| \ra 0$, it follows that $\la_n \ra \la$ (by the injectivity of the logarithm). 
Using the fact that $u\in \cc(\rd)$, we get 
\begin{align*}
   C_{y_n,\lambda_n}u(x) = \lambda_n^{-\frac{d-sp}{p}}u \left(\frac{x-y_n}{\lambda_n}\right) \ra \lambda^{-\frac{d-sp}{p}}u \left(\frac{x-y}{\lambda}\right) = C_{y,\la} u(x), 
\end{align*}
for every $x \in \rd$ as $n\to\infty$. Further, 
\begin{align*}
    \norm{ C_{y_n,\lambda_n}u }_{\wps} = \norm{u}_{\wps} = \norm{ C_{y,\lambda}u }_{\wps}.
\end{align*}
Therefore, applying Lemma \ref{convergence-BL}-(i), for all $u\in \cc(\rd)$ we get
\bee\no
\| C_{y_n,\lambda_n}u - C_{y,\lambda}u\|_{\wps}\to 0,\text{ as } n \to \infty. 
\eee  
Next, let $u\in \wps$ and $\varepsilon>0$ be given. By the density of $\cc(\rd)$ in $\wps$, there exists $v\in \cc(\rd)$ such that 
\begin{align*}
 \|u-v\|_{\wps} < \frac{\varepsilon}{3}.
\end{align*}
Hence there exists $n_1 \in \N$ such that for all $n \ge n_1$,
\begin{align*}
\|C_{y_n,\lambda_n}u - C_{y,\lambda}u\|_{\wps} &\leq \|C_{y_n,\lambda_n}u-C_{y_n,\lambda_n}v\|_{\wps} + \|C_{y_n,\lambda_n}v-C_{y,\lambda}v\|_{\wps} \no \\
&+ \|C_{y,\lambda}u - C_{y,\lambda}v\|_{\wps}\no\\
&= 2\|u-v\|_{\wps} + \|C_{y_n,\lambda_n}v-C_{y,\lambda}v\|_{\wps} < \varepsilon.
\end{align*}
Thus, $ C_{y_n,\lambda_n}u \rightarrow C_{y,\lambda}u$ for all $u\in \wps$, as required.
\end{proof}

In \cite[Lemma~3]{GiAd2014}, for $p= 2$, the authors established, if $C_{y_n, \la_n } u \rightharpoonup 0$ in $\mathcal{D}^{s,2}$ for every $u \in \mathcal{D}^{s,2}$, then $\left| \log(\lambda_n) \right| + |y_n| \to \infty$. In the same spirit, the following lemma shows that if a sequence in the group $G$ sends every element in $\wps$ to $0$ under the action $\mathcal{A}$, then the sequence must go to infinity with respect to the metric $d$ of $G$. 

\begin{proposition}\label{weak-bub}
Let $\{(y_n,\lambda_n)\} \subset G$ be such that
\begin{align}\label{gp-ele-infty}
\mathcal{A}(C_{y_n,\lambda_n}u,v)\to 0 \text{ for every }u,\,v\in\wps.
\end{align}
Then $\left| \log(\lambda_n) \right| + |y_n| \to \infty,$ as $n \ra \infty$.
Moreover, let $\{(a_n,\delta_n)\},\,\{(y_n,\lambda_n)\}\subset G$ be such that \begin{align}\label{gp-int-bub}
\mathcal{A}(C_{a_n,\delta_n}u,C_{y_n,\lambda_n}v)\to 0 \text{ for every }u,\,v\in\wps.
\end{align}
Then $\left| \log\left(\frac{\delta_n}{\lambda_n}\right) \right| + \left| \frac{a_n-y_n}{\lambda_n} \right| \to\infty,$ as $n \ra \infty$.
\end{proposition}

\begin{proof}
If \eqref{gp-ele-infty} fails, then up to subsequences we get $y_n\to y$ in $\rd$ and $\lambda_n\to\lambda>0$ as $n\to\infty$. Take $u\in \wps\setminus \{0\}$. Using Proposition \ref{Cylambda}, we have $C_{y_n, \lambda_n}(u) \ra C_{y, \lambda}(u)$ as $n \ra \infty$. By hypothesis, 
\bee\label{weak-conv-1}
\lim_{n \ra \infty} \mathcal{A}(C_{y_n,\lambda_n}u,w) = 0, \; \forall \,w\in\wps.
\eee
In particular, for $w=C_{y,\lambda}u\in \wps$, $\mathcal{A}(C_{y_n,\lambda_n}u, C_{y,\lambda}u) \to 0,\text{ as }n\to\infty.$
Since $\mathcal{A} \in C^1(\mathcal{W})$, we get
\begin{align*}
\norm{u}^p_{\wps} = \norm{C_{y,\lambda}u}_{\wps}^p = \mathcal{A}(C_{y,\lambda} u, C_{y,\lambda} u) = 0,
\end{align*}
which contradicts the fact that $u \neq 0$. Therefore, $\left| \log(\lambda_n) \right| + |y_n| \to \infty$ as $n \ra \infty$.  Further, using Lemma \ref{Asy-Bub-inter} for all $u,v\in\wps$, 
\begin{align*}
  \mathcal{A} \left( C_{\frac{a_n-y_n}{\lambda_n}, \frac{\delta_n}{\lambda_n}}u,v\right) = \mathcal{A} \left( C_{a_n,\delta_n}u,C_{y_n,\lambda_n}v \right) \to 0,\text{ as }n\to\infty.  
\end{align*}
Hence, in view of \eqref{gp-ele-infty}, we must have $\left| \log(\frac{\delta_n}{\lambda_n}) \right| + \left| \frac{a_n-y_n}{\lambda_n} \right| \to\infty$ as $n \ra \infty$.
\end{proof}

\begin{remark}\label{weak-bub-II}
Let $\{(y_n,\lambda_n)\} \subset G$ be a sequence such that for every $u \in \wps$, $C_{y_n,\lambda_n}u \rightharpoonup 0$ in $\wps$. Then $\mathcal{A}(C_{y_n,\lambda_n}u, v) \ra 0$ for every $u,v \in \wps$. As a consequence, using Proposition 
\ref{weak-bub}, we get
$\left| \log(\lambda_n) \right| + |y_n| \to \infty$, as $n \ra \infty$.
\end{remark}

Now we see the converse part of Proposition \ref{weak-bub}. In \cite[Lemma~3]{GiAd2014}, for $p= 2$, the authors also established that if $\left| \log(\lambda_n) \right| + |y_n| \to \infty$, then $C_{y_n, \la_n } u \rightharpoonup 0$ in $\mathcal{D}^{s,2}$ for every $u \in \mathcal{D}^{s,2}$. Their proof mainly uses stronger density results for $\mathcal{D}^{s,2}$, namely the density of $A\coloneqq \left\{ f \in \mathcal{S}(\rd): \mathcal{F}(f) \in \C_{c}^{\infty}(\rd \setminus \{ 0\}) \right\}$ in $\mathcal{D}^{s,2}$ (note that $\cc(\rd)$ is not contained in $A$). It follows that, if $u\in A$ then one has $(-\Delta)^su \in \mathcal{S}(\rd)$ (using the Fourier representation of $(-\Delta)^s$). However, for $p \neq 2$, we do not have a similar density result in $\wps$ and the Fourier representation of $(-\Delta_p)^s$. 
\begin{proposition}\label{weak-bub-converse}
 Let $s \in(0, \frac{1}{p'})$. Then for any sequence $\{(y_n,\lambda_n)\} \subset G$, we have 
\begin{align}\label{limit-2}
    \left| \log(\lambda_n) \right| + |y_n| \to \infty  \Longrightarrow \mathcal{A}(u,C_{y_n,\lambda_n}v)\to 0 \text{ and } \mathcal{A}(C_{y_n,\lambda_n}u, v)\to 0\text{ as }n\to\infty,
\end{align}
for all $u, v\in \wps$. 
\end{proposition}

\begin{proof}
First, we show that $\mathcal{A}(u,C_{y_n,\lambda_n}v)\to 0$ for all $u, v \in \wps$. 
Due to the density of $\cc(\rd)$ in $\wps$, it is enough to show $\mathcal{A}(u,C_{y_n,\lambda_n}v)\to 0$ for all $u, v \in \cc(\rd)$. If $\abs{y_n} \ra \infty$ and/or $\la_n \ra \infty$, then observe that $C_{y_n,\lambda_n}v \ra 0$ a.e. in $\rd$. Moreover, since $\{ C_{y_n,\lambda_n}v\}$ is bounded in $\wps$, from the reflexivity, up to a subsequence, $C_{y_n,\lambda_n}v \rightharpoonup w$ in $\wps$. The uniqueness of the limit yields $w=0$ a.e. in $\rd$.  Hence, using Lemma \ref{convergence-integrals}-(iii), we infer that $\mathcal{A}(u,C_{y_n,\lambda_n}v) \ra 0$ as $n \ra \infty$. Now we consider the case where $\la_n \ra 0$. In this case, using the change of variable $\bar{x} = x-y_n,\, \bar{y} = y-y_n$ we write 
\begin{align*}
    &\mathcal{A}(u, C_{y_n,\lambda_n}v) = \la_n^{-\frac{d-sp}{p}} \iint\limits_{\rd \times \rd} \frac{\abs{u(x) - u(y)}^{p-2}(u(x) - u(y)) \left( v(\frac{x-y_n}{\la_n}) - v(\frac{y-y_n}{\la_n}) \right)}{\abs{x-y}^{d+sp}} \dxy \\
    & = \la_n^{-\frac{d-sp}{p}} \iint\limits_{\rd \times \rd} \frac{\abs{u(x+y_n) - u(y+y_n)}^{p-2}(u(x+y_n) - u(y+y_n)) \left( v(\frac{x}{\la_n}) - v(\frac{y}{\la_n}) \right)}{\abs{x-y}^{d+sp}} \dxy \\
    &= I_1 + I_2, 
\end{align*}
where 
\begin{align*}
    &I_1\coloneqq \la_n^{-\frac{d-sp}{p}} \iint\limits_{\rd \times \rd} \frac{\abs{u(x+y_n) - u(y+y_n)}^{p-2}(u(x+y_n) - u(y+y_n)) v(\frac{x}{\la_n}) }{\abs{x-y}^{d+sp}} \dxy, \\
    &I_2\coloneqq \la_n^{-\frac{d-sp}{p}} \iint\limits_{\rd \times \rd} \frac{\abs{u(y+y_n) - u(x+y_n)}^{p-2}(u(y+y_n) - u(x+y_n)) v(\frac{y}{\la_n}) }{\abs{x-y}^{d+sp}} \dxy. 
\end{align*}
We show that $I_1 \ra 0$ as $n \ra \infty$. By interchanging the role of $x$ and $y$, $I_2 \ra 0$ also follows. Since $u, v \in \cc(\rd)$, there exists $R>0$ such that $\text{supp}(u),\, \text{supp}(v) \subset B(0, R)$. We split $I_1$ as follows: 
\begin{align*}
    I_1 & \coloneqq \la_n^{-\frac{d-sp}{p}} \left(  \iint\limits_{B(0,R) \times B(0,R)} + \iint\limits_{B(0,R) \times B(0,R)^c} + \iint\limits_{B(0,R)^c \times \rd} \right) \\
    &\frac{\abs{u(x+y_n) - u(y+y_n)}^{p-2}(u(x+y_n) - u(y+y_n)) v(\frac{x}{\la_n}) }{\abs{x-y}^{d+sp}} \dxy \coloneqq I_{1,1} + I_{1,2} + I_{1,3}. 
\end{align*}
First, we estimate $I_{1,3}$. Since $\abs{x} \ge R$ and $\la_n \ra 0$, for large $n \in \N$, we have 
\begin{align*}
    \left| \frac{x}{\la_n} \right| \ge \frac{R}{\la_n} > R,
\end{align*}
and hence $v(\frac{x}{\la_n}) = 0$. Therefore, $I_{1,3} = 0$ for large $n \in \N$. 

\noi \textbf{Estimation of $I_{1,1}$:} Now we estimate $I_{1,1}$. Set $h=y-x$. Using $\abs{h} \le \abs{y} + \abs{x} \le 2R$, we get 
\begin{align}\label{estimate-1}
    \int_{B(0, R)} \frac{\abs{u(x+y_n) - u(y+y_n)}^{p-1}}{\abs{x-y}^{d+sp}} \dy \le \int_{B(0,2R)} \frac{\abs{u(x+y_n) - u(x+h+y_n)}^{p-1}}{\abs{h}^{d+sp}} \dh.
\end{align}
Since $u \in \cc(\rd)$, there exists $C>0$ such that $\abs{u(x+h+y_n)-u(x+y_n)} \le C \abs{h}$. Hence the integral in \eqref{estimate-1} can be further estimated by
\begin{align*}
    \int_{B(0, 2R)} \frac{\dh}{\abs{h}^{d+sp-p+1}} = C(R) < \infty \text{ if } s \in (0, \frac{1}{p'}). 
\end{align*}
Therefore, if $s \in (0, \frac{1}{p'})$, then we have
\begin{align}\label{I-11}
    I_{1,1} \lesssim_{R} \la_n^{-\frac{d-sp}{p}} \int_{B(0, R)} \left| v(\frac{x}{\la_n}) \right| \dx. 
\end{align}

Next, we prove the following claim:

\noi \textbf{Claim:} For any finite $R_1>0$ and $\kappa \in [1, p^*_s)$, 
\begin{align}\label{claim-1}
    \la_n^{-\frac{d-sp}{p}} \left( \int_{B(0, R_1)} \left| v(\frac{x}{\la_n}) \right|^{\kappa} \dx \right)^{\frac{1}{\kappa}} \ra 0, \text{ as } n \ra \infty.
\end{align}

\noi To prove the claim, we use the change of variable $x=\la_n y$ to get
\begin{align*}
    \la_n^{-\frac{d-sp}{p}} \left( \int_{B(0, R_1)} \left| v(\frac{x}{\la_n}) \right|^{\kappa} \dx \right)^{\frac{1}{\kappa}} & = \la_n ^{\left( \frac{d}{\kappa} - \frac{d-sp}{p} \right)} \left( \int_{B(0, \frac{R_1}{\la_n} )} \left| v(y) \right|^{\kappa} \dy \right)^{\frac{1}{\kappa}} \\
    & \le \la_n ^{\left( \frac{d}{\kappa} - \frac{d-sp}{p} \right)} \left( \int_{\rd} \left| v(y) \right|^{\kappa} \dy \right)^{\frac{1}{\kappa}},
\end{align*}
where using the fact that $\kappa <p^*_s \Longleftrightarrow \frac{d}{\kappa}>\frac{d-sp}{p}$,  
\begin{align*}
    \la_n^{\left( \frac{d}{\kappa} - \frac{d-sp}{p} \right)} \ra 0, \text{ as } n \ra \infty,  
\end{align*}
and further $\int_{\rd} \left| v(y) \right|^{\kappa} \dy < \infty$ as $v \in \C_c^{\infty}(\rd)$. This proves the claim. 

Using \eqref{I-11} and \eqref{claim-1} we conclude $I_{1,1} \ra 0$ as $n \ra \infty$.

\noi \textbf{Estimation of $I_{1,2}$:} Next, to estimate $I_{1,2}$, we write 
\begin{align}\label{I-12}
    I_{1,2} \le \la_n^{-\frac{d-sp}{p}} \left( \iint\limits_{B(0,R) \times A} + \iint\limits_{B(0,R) \times B(0,R)^c \setminus A}\right) \frac{\abs{u(x+y_n) - u(y+y_n)}^{p-1} \left| v(\frac{x}{\la_n}) \right| }{\abs{x-y}^{d+sp}} \dxy,
\end{align}
where $A\coloneqq  \{ y \in B(0, R)^c : \abs{x-y} \le 1 \}$. For $y \in A,\, x \in B(0,R)$, $\abs{y} \le \abs{y-x} + \abs{x} \le R+1$. Hence, the first integral can be estimated as 
\begin{align*}
    \la_n^{-\frac{d-sp}{p}} \iint\limits_{B(0,R) \times B(0,R+1)} \frac{\abs{u(x+y_n) - u(y+y_n)}^{p-1} \left| v(\frac{x}{\la_n}) \right| }{\abs{x-y}^{d+sp}} \dxy \ra 0,
\end{align*}
as $n \ra \infty$, where the convergence follows using the same arguments as in the estimation of $I_{1,1}$. Finally, using the change of variable $y=x-h$ and Fubini's theorem, we estimate the second integral of \eqref{I-12} as 
\begin{align*}
  &\la_n^{-\frac{d-sp}{p}}  \int_{\{\abs{h}>1\}} \frac{\dh}{\abs{h}^{d+sp}} \left( \int_{B(0, R)} \chi_{B(x,R)^c}  \abs{u(x+y_n) - u(x-h+y_n)}^{p-1} \left| v(\frac{x}{\la_n}) \right| \dx  \right) \\
  &\le \norm{u}_{L^{\infty}(\rd)}^{p-1} \left( \int_{\{\abs{h}>1\}} \frac{\dh}{\abs{h}^{d+sp}} \right)  \left( \la_n^{-\frac{d-sp}{p}} \int_{B(0, R)} \left| v(\frac{x}{\la_n}) \right| \dx  \right) \ra 0, \text{ as } n \ra \infty,
\end{align*}
using \eqref{claim-1} and the fact that $\int_{\{\abs{h}>1\}} \frac{\dh}{\abs{h}^{d+sp}} < \infty$. Therefore, $I_{1,2} \ra 0$.

Hence $\mathcal{A}(u, C_{y_n,\lambda_n}v) \ra 0$ as $n \ra \infty$. Now, using Lemma \ref{Asy-Bub-inter},
\begin{align*}
    \mathcal{A}(C_{y_n,\lambda_n}u,v) = \mathcal{A}(u,C_{-\lambda_n^{-1}y_n,\lambda_n^{-1}}v), \; \forall \, (u,v) \in \mathcal{W}.
\end{align*}  
Set $\delta_n = \frac{1}{\la_n}$ and $x_n = - \frac{y_n}{\la_n}$. Now observe that 
\begin{align*}
    \abs{\log(\delta_n)} + \abs{x_n} \ra \infty \Longleftrightarrow \abs{\log(\la_n)} + \abs{y_n} \ra \infty.
\end{align*}
Therefore, we conclude that \eqref{limit-2} holds.
\end{proof}

In view of Lemma \ref{Asy-Bub-inter} and Proposition \ref{weak-bub-converse}, we have the following. 

\begin{proposition}\label{weak-bub-converse-II}
Let $s \in(0, \frac{1}{p'})$. For any sequence $\{(a_n,\delta_n)\},\,\{(y_n,\lambda_n)\}\subset G$ we have
\begin{align*}
\left| \log\left(\frac{\delta_n}{\lambda_n}\right) \right| + \left| \frac{a_n-y_n}{\lambda_n} \right| \to\infty \Longrightarrow \mathcal{A}(C_{a_n,\delta_n}u,C_{y_n,\lambda_n}v)\to 0,
\end{align*}
as $n \ra \infty$ and for all $u,\,v\in\wps$. 
\end{proposition}

Now, we consider the following map
\begin{align}\label{action-2}
    \mathcal{A}_1(u,v) \coloneqq \iint\limits_{\rd \times \rd} \frac{\abs{u(x) - u(y)}^{p-1} \abs{v(x) - v(y)}}{\abs{x-y}^{d+sp}} \dxy, \; \forall \, (u,v) \in \mathcal{W}.
\end{align}
Using a similar set of arguments as in Proposition \ref{weak-bub-converse-II} (with necessary modifications):

\begin{lemma}\label{weak-bub-converse-III}
 Let $s \in(0, \frac{1}{p'})$. For any sequence $\{(a_n,\delta_n)\},\,\{(y_n,\lambda_n)\}\subset G$ we have
\begin{align*}
\left| \log\left(\frac{\delta_n}{\lambda_n}\right) \right| + \left| \frac{a_n-y_n}{\lambda_n} \right| \to\infty \Longrightarrow \mathcal{A}_1(C_{a_n,\delta_n}u,C_{y_n,\lambda_n}v)\to 0,
\end{align*}
as $n \ra \infty$ and for all $u,\,v\in\wps$.
\end{lemma}

\subsection{Refined embeddings of \texorpdfstring{$\wps$}{} in the Morrey spaces}

Now, we discuss some embeddings of $\wps$ in some Morrey spaces, which are finer than the classical embeddings $\wps \hookrightarrow L^{p^*_s}(\rd)$. First, we recall the definition of Morrey space introduced by Morrey as a refinement of Lebesgue space. Let $\mathcal{M}(\rd)$ be the set of all measurable functions defined on $\rd$.
\begin{definition}[Morrey space]\label{Morreyspaces} For $r \in [1, \infty)$ and $\gamma \in [0, d]$, the homogeneous Morrey space $\mathcal{L}^{r,\gamma}(\rd)$ is defined as 
  \begin{equation*}
    \mathcal{L}^{r,\gamma}(\mathbb{R}^d) \coloneqq \left\{u \in \mathcal{M}(\rd): \sup_{x\in\rd,\, {R>0}}R^{\gamma}\fint_{B(x,R)} |u(y)|^r \dy < \infty \right\},
\end{equation*}
with norm 
\bee\no
\norm{u}_{\mathcal{L}^{r,\gamma}(\rd)}^r \coloneqq \sup_{x\in\rd,\, {R>0}}R^{\gamma}\fint_{B(x,R)} |u(y)|^r \dy.
\eee
\end{definition}
Note, that if $\gamma =d$ then $\mathcal{L}^{r,d}(\rd)=L^{r}(\rd)$ for $r\geq 1$ and whenever $\gamma=0$, $\mathcal{L}^{r,0}(\rd)=L^{\infty}(\rd)$. Observe that whenever $r\in [1,p_s^*]$, 
\bee\no
\|C_{y,\lambda}u\|_{\mathcal{L}^{r,\frac{d-sp}{p}r}(\rd)} = \|u\|_{\mathcal{L}^{r,\frac{d-sp}{p}r}(\rd)},\text{ for }y\in\rd \text{ and }\lambda>0,  
\eee

\begin{remark}
    Using H\"older's inequality and definition of Morrey Space, the following chain of continuous embeddings hold:
\begin{align*}
        \wps \hookrightarrow L^{p_s^*}(\rd) \hookrightarrow \mathcal{L}^{r,\frac{r(d-sp)}{p}}(\rd) \hookrightarrow \mathcal{L}^{p,d-sp}(\rd) \hookrightarrow \mathcal{L}^{1,\frac{d-sp}{p}}(\rd).
    \end{align*}
  where $r \in [p, p^*_s]$.
\end{remark}

In \cite[Lemma 3.1]{Zhang-2021}, Zhang proved that 
\begin{align}\label{finer-1}
   \norm{u}_{L^{p_s^*}(\rd)} \lesssim_{d,p,s} \norm{u}_{\wps}^{\frac{p}{p_s^*}} \norm{u}_{\mathcal{L}^{p,d-sp}(\rd)}^{1-\frac{p}{p_s^*}}, \; \forall\, u \in \wps. 
\end{align}
Applying $L^{p_s^*}(\rd) \hookrightarrow \mathcal{L}^{p,d-sp}(\rd)$ one clearly see that \eqref{finer-1} yields $\wps(\rd) \hookrightarrow L^{p^*_s}(\rd)$.
\begin{proposition}\label{interpolation(product)}
Let $\frac{p}{p^*_s} \le \theta <1$ and $p \le r \le p^*_s$. Then 
\bee\label{Inter-morrey-sobolev}
\norm{u}_{L^{p^*_s}(\rd) } \lesssim_{d,p,s} \norm{u}_{\wps}^{\theta} \norm{u}_{\mathcal{L}^{r, \frac{r(d-sp)}{p}}(\rd)}^{1-\theta}, \; \forall \, u \in \wps.
\eee
Moreover, 
\begin{align}\label{Inter-morrey-sobolev-2}
\norm{(u,v)}_{L^{p^*_s}(\rd) \times L^{p^*_s}(\rd) } \lesssim_{d,p,s} \norm{(u,v)}_{\W}^{\theta} \norm{(u,v)}_{\mathcal{L}^{r, \frac{r(d-sp)}{p}}(\rd) \times \mathcal{L}^{r, \frac{r(d-sp)}{p}}(\rd)}^{1-\theta}, \; \forall \, (u,v) \in \W.
\end{align}
\end{proposition}

\begin{proof}
  Using the embedding $\mathcal{L}^{r,\frac{r(d-sp)}{p}}(\rd) \hookrightarrow \mathcal{L}^{p,d-sp}(\rd)$ and \eqref{finer-1}, we get 
     \begin{align*}
    \|u\|_{L^{p_s^*}(\rd)} \lesssim_{d,p,s} \norm{u}_{\wps}^{\frac{p}{p_s^*}} \norm{u}_{\mathcal{L}^{r,\frac{r(d-sp)}{p}}(\rd)}^{\frac{p}{p_s^*}}, \; \forall \, u \in \wps.
     \end{align*}
    Now for $\delta\in(0,1-\frac{p}{p_s^*})$, using $\wps \hookrightarrow \mathcal{L}^{r,\frac{r(d-sp)}{p}}(\rd)$, 
    \begin{align*}
        \|u\|_{L^{p_s^*}(\rd)} & \lesssim_{d,p,s}  \norm{u}_{\wps}^{\frac{p}{p_s^*}} \|u\|_{\mathcal{L}^{r,\frac{r(d-sp)}{p}}(\rd)}^{\delta} \|u\|_{\mathcal{L}^{r,\frac{r(d-sp)}{p}}(\rd)}^{1-\frac{p}{p_s^*}-\delta} \\ &\lesssim_{d,p,s} \norm{u}_{\wps}^{\frac{p}{p_s^*}+\delta} \|u\|_{\mathcal{L}^{r,\frac{r(d-sp)}{p}}(\rd)}^{1-\frac{p}{p_s^*}-\delta}, \; \forall \, u \in \wps. 
    \end{align*}
    Taking $\theta = \frac{p}{p_s^*} +\delta$, we see that \eqref{Inter-morrey-sobolev} holds for all $\theta \in [\frac{p}{p_s^*},1)$. Further, \eqref{Inter-morrey-sobolev-2} holds using $$\norm{(u,v)}_{L^{p_s^*}(\rd) \times L^{p_s^*}(\rd)} \le \norm{u}_{L^{p_s^*}(\rd)} + \norm{v}_{L^{p_s^*}(\rd)}, \; \forall\, u,v \in L^{p_s^*}(\rd),$$ and \eqref{Inter-morrey-sobolev}. 
  \end{proof}

\subsection{The Palais-Smale decomposition}
This section studies the Palais-Smale decomposition of $I_{f,g}$. 

\noi \textbf{Proof of Theorem \ref{PS-decomposition}:}
Since $\{ (u_n,v_n) \} \subset \W$  is a (PS) sequence for $I_{f,g}$ at level $\eta$,  
\begin{align}\label{PSD-1}
    I_{f,g}(u_n, v_n) - \frac{1}{p^*_s}  \prescript{}{\W'}{\langle} I'_{f,g}(u_n, v_n),(u_n, v_n){\rangle}_{\W} \le \eta + o_n(1) + o_n(1) \norm{(u_n, v_n)}_{\W}.
\end{align}
Now 
\begin{align*}
    & \text{ L.H.S. of \eqref{PSD-1} } \ge \left( \frac{1}{p} - \frac{1}{p^*_s} \right) \norm{(u_n, v_n)}_{\W}^p \\
    & \quad \quad - \left( 1 - \frac{1}{p^*_s} \right) \left( \norm{f}_{(\wps)'} \norm{u_n}_{\wps} + \norm{g}_{(\wps)'} \norm{v_n}_{\wps} \right) \\
    & \ge \left( \frac{1}{p} - \frac{1}{p^*_s} \right) \norm{(u_n, v_n)}_{\W}^p -  \left( 1 - \frac{1}{p^*_s} \right) \left( \norm{f}_{(\wps)'} + \norm{g}_{(\wps)'} \right) \norm{(u_n, v_n)}_{\W}.
\end{align*}
In view of R.H.S. of \eqref{PSD-1}, $\{ (u_n, v_n ) \}$ is a bounded sequence on $\W$. By the reflexivity of $\W$, let $\{ (u_n, v_n) \}$ weakly converge to $(\tilde{u},\tilde{v})$ in $\W$ (up to a subsequence). Since $I'_{f, g}(u_n,v_n) \ra 0$ in $\mathcal{W}'$, for every $(\phi, \psi) \in \W$ we have
\begin{align*}
    \mathcal{A}(u_n , \phi) + \mathcal{A}(v_n, \psi) & = \frac{\al}{p^*_s} \int_{\rd} \abs{u_n}^{\al-2}u_n\abs{v_n}^{\beta} \phi \dx + \frac{\be}{p^*_s} \int_{\rd} \abs{v_n}^{\be-2}v_n\abs{u_n}^{\al} \psi \dx \no \\
    &+\prescript{}{(\wps)'}{\langle}f,\phi{\rangle}_{\wps}+\prescript{}{(\wps)'}{\langle}g,\psi{\rangle}_{\wps} + o_n(1). 
\end{align*}
Taking the limit as $n \ra \infty$ in the above identity and using Lemma \ref{convergence-integrals}, we see that $(\tilde{u},\tilde{v})\in \W$ satisfies \eqref{weak1} weakly. We divide the rest of the proof into several steps. 

\noi \textbf{Step 1:} In this step, we claim that $\{ u_n - \tilde{u}, v_n - \tilde{v}\}$ is a (PS) sequence for $I_{0,0}$ at level $\eta -I_{f,g}(\tilde{u}, \tilde{v})$. Set $\tilde{u}_n = u_n - \tilde{u}, \tilde{v}_n = v_n - \tilde{v}$. Using Lemma \ref{convergence2}, Lemma \ref{convergence-BL}, and $(\tilde{u}_n,\tilde{v}_n) \rightharpoonup (0,0)$, we get
\begin{align*}
    I_{0,0}(\tilde{u}_n, \tilde{v}_n) & = \frac{1}{p} \norm{(\tilde{u}_n, \tilde{v}_n)}_{\W}^p - \frac{1}{p^*_s} \int_{\rd}\abs{\tilde{u}_n}^{\alpha} \abs{\tilde{v}_n}^{\beta} \\
    & = \frac{1}{p} \left( \norm{u_n}_{\wps}^p - \norm{\tilde{u}}_{\wps}^p \right) + \frac{1}{p} \left( \norm{v_n}_{\wps}^p - \norm{\tilde{v}}_{\wps}^p \right) - \frac{1}{p^*_s}  \int_{\rd} \left( \abs{u_n}^{\al} \abs{v_n}^{\be} -  \abs{\tilde{u}}^{\al} \abs{\tilde{v}}^{\be}\right) \\
    & -\prescript{}{(\wps)'}{\langle}f,u_n{\rangle}_{\wps}-\prescript{}{(\wps)'}{\langle}g,v_n{\rangle}_{\wps} + \prescript{}{(\wps)'}{\langle}f, \tilde{u}{\rangle}_{\wps}+ \prescript{}{(\wps)'}{\langle}g, \tilde{v}{\rangle}_{\wps} + o_n(1) \\
    & = I_{f,g}(u_n, v_n) - I_{f,g}(\tilde{u}, \tilde{v}) + o_n(1).
\end{align*}
Hence $I_{0,0}(\tilde{u}_n, \tilde{v}_n) \ra \eta - I_{f,g}(\tilde{u}, \tilde{v})$ as $n \ra \infty$. Further, for $(\phi, \psi) \in \W$, using Lemma \ref{convergence-integrals} we have 
\begin{align*}
    &\prescript{}{\W'}{\langle} I_{0,0}'(\tilde{u}_n, \tilde{v}_n), (\phi, \psi) {\rangle}_{\W} \\
    &= \mathcal{A}(\tilde{u}_n , \phi) + \mathcal{A}(\tilde{v}_n,\psi) - \frac{\al}{p^*_s} \int_{\rd} \abs{\tilde{u}_n}^{\al -2} \tilde{u}_n \abs{\tilde{v}_n}^{\be} \phi - \frac{\be}{p^*_s} \int_{\rd} \abs{\tilde{v}_n}^{\be-2} \tilde{v}_n \abs{\tilde{u}_n}^{\al} \psi \rightarrow 0, \text{ as } n \ra \infty.
\end{align*}
Thus, the claim holds. 

\noi \textbf{Step 2:} Suppose $(u_n,v_n) \ra (\tilde{u}, \tilde{v})$ as $n\to\infty$ in $\W$. Then using the continuity of $I_{f,g}$, we get $\eta = I_{f,g}(\tilde{u}, \tilde{v})$, and Theorem \ref{PS-decomposition} holds for $k=0$. So, we assume that  $(u_n,v_n) \not\ra (\tilde{u}, \tilde{v})$ in $\W$. In view of Step 1, $\prescript{}{\W'}{\langle} I_{0,0}'(\tilde{u}_n, \tilde{v}_n), (\tilde{u}_n, \tilde{v}_n) {\rangle}_{\W} \ra 0$, which implies 
\begin{align}\label{del-0}
    0 < c \le  \norm{(\tilde{u}_n, \tilde{v}_n)}^p_{\W} = \int_{\rd} \abs{\tilde{u}_n}^{\alpha} \abs{\tilde{v}_n}^{\beta} + o_n(1).
\end{align}
In this step, we construct a sequence $(\hat{u}_n, \hat{v}_n)$ from $(\tilde{u}_n, \tilde{v}_n)$ in such a way that the $\wps$-norm is preserved and $(\hat{u}_n, \hat{v}_n)$ weakly goes to a non-zero limit $(\hat{u}, \hat{v}) \in \W$. From \eqref{del-0} we have 
\begin{align*}
    0 < c \le  \norm{(\tilde{u}_n, \tilde{v}_n)}^p_{\W} = \int_{\rd} \abs{\tilde{u}_n}^{\alpha} \abs{\tilde{v}_n}^{\beta} \le \int_{\rd} \abs{\tilde{u}_n}^{p^*_s} + \int_{\rd} \abs{\tilde{v}_n}^{p^*_s} = \norm{(\tilde{u}_n, \tilde{v}_n)}_{L^{p^*_s}(\rd) \times L^{p^*_s}(\rd)}^{p^*_s}.
\end{align*}
Therefore, $(\tilde{u}_n, \tilde{v}_n) \not \ra 0$ in $L^{p^*_s}(\rd) \times L^{p^*_s}(\rd)$, and hence there exists $\delta>0$ such that 
\begin{align*}
    \inf_{n \in \mathbb{N}} \norm{(\tilde{u}_n, \tilde{v}_n)}_{L^{p^*_s}(\rd) \times L^{p^*_s}(\rd)}^{p^*_s} \ge \delta.
\end{align*}
Hence the interpolation given in \eqref{Inter-morrey-sobolev-2} yields 
\begin{align*}
    \delta \le C \norm{(\tilde{u}_n, \tilde{v}_n)}_{\W}^{\theta} \norm{(\tilde{u}_n, \tilde{v}_n)}_{\mathcal{L}^{p,d-sp}(\rd) \times \mathcal{L}^{p,d-sp}(\rd)}^{1- \theta} \le C_1 \norm{(\tilde{u}_n, \tilde{v}_n)}_{\mathcal{L}^{p,d-sp}(\rd) \times \mathcal{L}^{p,d-sp}(\rd)}^{1- \theta},
\end{align*}
which implies $\norm{(\tilde{u}_n, \tilde{v}_n)}_{\mathcal{L}^{p,d-sp}(\rd) \times \mathcal{L}^{p,d-sp}(\rd)} \ge C_2$ for some $C_2>0$, i.e., 
\begin{align*}
    \sup_{x \in \rd,\, R>0} R^{d-sp} \fint_{B(x,R)} \left( |\tilde{u}_n(y)|^p + |\tilde{v}_n(y)|^p \right) \dy \ge C_2.  
\end{align*}
Hence for every $n \in \N$, there exist $x_n \in \rd$ and $r_n \in (0, \infty)$ such that
\begin{align}\label{PSD-3}
    r_n^{d-sp} \fint_{B(x_n, r_n)} \left( |\tilde{u}_n(y)|^p + |\tilde{v}_n(y)|^p \right) \dy \ge \norm{(\tilde{u}_n, \tilde{v}_n)}_{\mathcal{L}^{p,d-sp}(\rd) \times \mathcal{L}^{p,d-sp}(\rd)}^p - \frac{C_2}{2n} > C_3,
\end{align}
for some $C_3>0$. Set 
\begin{align*}
    \hat{u}_n(z) = r_n^{\frac{d-sp}{p}} \tilde{u}_n(r_n z + x_n) \text{ and } \hat{v}_n(z) = r_n^{\frac{d-sp}{p}} \tilde{v}_n(r_n z + x_n), \text{ for } z \in \rd.
\end{align*}
Let $(\hat{u}_n, \hat{v}_n) \rightharpoonup (\hat{u}, \hat{v})$ in $\W$. By the change of variable and using the compact embedding of $\W \hookrightarrow L_{loc}^p(\rd) \times L_{loc}^p(\rd)$, we get 
\begin{align*}
    \omega_N C_3 \leq r_n^{-sp}  \int_{B(x_n, r_n)} \left( |\tilde{u}_n(y)|^p + |\tilde{v}_n(y)|^p \right) \dy & = \int_{B(0,1)} \left( |\hat{u}_n(y)|^p + |\hat{v}_n(y)|^p \right) \dy \\
    & \longrightarrow \int_{B(0,1)} \left( |\hat{u}(y)|^p + |\hat{v}(y)|^p \right) \dy, \text{ as } n \ra \infty.
\end{align*}
Therefore, $\hat{u}, \hat{v} \neq 0$. Next, we prove the following claim:

\noi \textbf{Claim:} If  $\{x_n\}$ is bounded, then it holds
\begin{align}\label{r_n limit}
    \lim_{n \ra \infty} r_n \in \{0,\infty\}.
\end{align}
Suppose $r_n \ra r_0 \in (0,\infty)$. Then $r_n^{-sp} \ra r_0^{-sp}$. 
Set 
\begin{align*}
    U_n(z) = r_n^{-sp} \chi_{B(x_n, r_n)}(z) \abs{\tilde{u}_n(z)}^p \text{ and } V_n(z) = r_n^{-sp} \chi_{B(x_n, r_n)}(z) \abs{\tilde{v}_n(z)}^p, \text{ for } z \in \rd.
\end{align*}
Now $(\tilde{u}_n, \tilde{v}_n) \rightharpoonup (0,0)$ in $\W$ and the compact embedding $\wps \hookrightarrow L_{loc}^p(\rd)$ yield $(\tilde{u}_n, \tilde{v}_n) \ra (0,0)$ in $L_{loc}^p(\rd) \times L_{loc}^p(\rd)$, and up to a subsequence $\tilde{u}_n(x)\to 0$ and $\tilde{v}_n(x)\to 0$ for a.e. $x\in\rd$. Suppose $x_n \ra x_0$ for some $x_0 \in \rd$. Then there exists $R$ large enough and $n_1 \in \N$ such that $B(x_n, r_n) \subset B(x_0, R)$ for $n \ge n_1$. Clearly, for $n\geq n_1,\,\text{supp}(U_n) \subset B(x_0, R),\, \text{supp}(V_n) \subset B(x_0, R)$ and $(U_n,V_n) \ra (0,0)$ in $L^1(B(x_0, R)) \times L^1(B(x_0, R))$. Hence 
\begin{align*}
    r_n^{-sp} \int_{B(x_n, r_n)} \left( \abs{\tilde{u}_n(y)}^p + \abs{\tilde{v}_n(y)}^p \right)\dy & = \int_{\rd} \left( \abs{U_n(y)} + \abs{V_n(y)} \right)\dy \\
    & = \int_{B(x_0, R)} \left( \abs{U_n(y)} + \abs{V_n(y)} \right)\dy \rightarrow 0,
\end{align*}
as $n \ra \infty$, which again contradicts \eqref{PSD-3}. 

Thus, in conclusion, up to a subsequence, either $x_n\to x_0$ and $\lim_{n\to\infty}r_n \in \{0, \infty\}$, or $|x_n|\to\infty$.

\noi \textbf{Step 3:} In this step, we claim that $(\hat{u}, \hat{v})$ weakly solves the following homogeneous system: 
\begin{equation*}
\left\{\begin{aligned}
&(-\Delta_p)^s \hat{u} = \frac{\al}{p_{s}^{*}} |\hat{u}|^{\alpha-2} \hat{u} |\hat{v}|^{\beta}\;\;\text{in}\;\rd,\\
&(-\Delta_p)^s \hat{v} = \frac{\beta}{p^{*}_{s}} |\hat{v}|^{\beta-2} \hat{v} |\hat{u}|^{\alpha}\;\;\text{in}\;\rd.
\end{aligned}
\right.
\end{equation*}
Take $\phi, \psi \in \wps$.  From Step 2, since $\hat{u}_n \rightharpoonup \hat{u}$ and $\hat{v}_n \rightharpoonup \hat{v}$ in $\wps$, using Lemma \ref{convergence-integrals}-(ii) we get 
\begin{align}\label{PSD_3}
   \lim_{n\to\infty} \mathcal{A}(\hat{u}_n, \phi) = \mathcal{A}(\hat{u}, \phi) \text{ and }  \lim_{n\to\infty}\mathcal{A}(\hat{v}_n, \psi) = \mathcal{A}(\hat{v}, \psi). 
\end{align}
For $n \in \mathbb{N}$, we set 
\begin{align*}
    \phi_n(z) = r_n^{- \frac{d-sp}{p}} \phi \left(\frac{z -  x_n}{r_n} \right) \text{ and } \psi_n(z) = r_n^{- \frac{d-sp}{p}} \psi \left(\frac{z -  x_n}{r_n} \right), \text{ for } z \in \rd.
\end{align*}
Note that $\norm{\phi_n}_{\wps} = \norm{\phi}_{\wps}$ and $\norm{\psi_n}_{\wps} = \norm{\psi}_{\wps}$.
Next, using the change of variable $\overline{x}_n = r_n x + x_n,\, \overline{y}_n = r_n y + y_n$, 
\begin{align}\label{rngoesto0}
    & \mathcal{A}(\hat{u}_n, \phi) = r_n^{\frac{d-sp}{p'}} \iint\limits_{\rd \times \rd} \frac{\abs{\tilde{u}_n(r_n x + x_n) - \tilde{u}_n(r_n y + x_n)}^{p-2} (\tilde{u}_n(r_n x + x_n) - \tilde{u}_n(r_n y + x_n))}{\abs{x-y}^{d+sp}} \no \\
    &\quad\quad\quad\quad\quad\quad\quad\quad(\phi(x) - \phi(y)) \dx\dy \no \\
    & = r_n^{-\frac{d-sp}{p}} \iint\limits_{\rd \times \rd} \frac{\abs{\tilde{u}_n(r_n x + x_n) - \tilde{u}_n(r_n y + x_n)}^{p-2} (\tilde{u}_n(r_n x + x_n) - \tilde{u}_n(r_n y + x_n))(\phi(x) - \phi(y))}{\abs{(r_n x + x_n)-(r_n y + x_n)}^{d+sp}}  \dx\dy  \no \\
    & = r_n^{- \frac{d-sp}{p}} \iint\limits_{\rd \times \rd} \frac{\abs{\tilde{u}_n(\overline{x}_n) - \tilde{u}_n(\overline{y}_n)}^{p-2} (\tilde{u}_n(\overline{x}_n) - \tilde{u}_n(\overline{y}_n)) \left(\phi \left(\frac{\overline{x}_n -  x_n}{r_n} \right) - \phi \left(\frac{\overline{y}_n -  x_n}{r_n} \right) \right)}{\abs{\overline{x}_n - \overline{y}_n}^{d+sp}} \d \overline{x}_n \d \overline{y}_n \no \\
    &= \mathcal{A}(\tilde{u}_n, \phi_n),
\end{align}
and proceeding similarly as in \eqref{rngoesto0} we also get
$\mathcal{A}(\hat{v}_n, \psi) = \mathcal{A}(\tilde{v}_n, \psi_n)$. 
Now using the fact that $\prescript{}{\W'}{\langle} I_{0,0}'(\tilde{u}_n, \tilde{v}_n), (\phi_n, \psi_n) {\rangle}_{\W} \ra 0$, we obtain
\begin{align}\label{PSD-4}
    &\mathcal{A}(\hat{u}_n, \phi) + \mathcal{A}(\hat{v}_n, \psi) = \mathcal{A}(\tilde{u}_n, \phi_n) + \mathcal{A}(\tilde{v}_n, \psi_n) \no \\
    &= \int_{\rd} \abs{\tilde{u}_n( \overline{x}_n)}^{\al-2}\tilde{u}_n(\overline{x}_n) \abs{\tilde{v}_n(\overline{x}_n)}^{\beta} \phi_n(\overline{x}_n) \d \overline{x}_n  + \int_{\rd} \abs{\tilde{v}_n(\overline{x}_n)}^{\beta-2} \tilde{v}_n(\overline{x}_n) \abs{\tilde{u}_n( \overline{x}_n)}^{\al} \psi_n(\overline{x}_n) \d \overline{x}_n + o_n(1) \no \\
    & = \int_{\rd} \abs{\hat{u}_n( x)}^{\al-2}\hat{u}_n(x) \abs{\hat{v}_n(x)}^{\beta} \phi(x) \dx + \int_{\rd} \abs{\hat{v}_n( x)}^{\be-2}\hat{v}_n(x) \abs{\hat{u}_n(x)}^{\al} \psi(x) \dx + o_n(1),
\end{align}
where the last identity holds from the following arguments 
\begin{align*}
    &\int_{\rd} \abs{\hat{u}_n( x)}^{\al-2}\hat{u}_n(x) \abs{\hat{v}_n(x)}^{\beta} \phi(x) \dx \\
    &= r_n^{\frac{d-sp}{p}(p^*_s -1)} \int_{\rd} \abs{\tilde{u}_n(r_n x + x_n)}^{\al -2} \tilde{u}_n(r_n x + x_n) \abs{\tilde{v}_n(r_n x + x_n)}^{\beta} \phi(x) \dx \\
    & = \frac{r_n^{\frac{dp-d+sp}{p}}}{r_n^d} \int_{\rd} \abs{\tilde{u}_n(\overline{x}_n)}^{\al -2} \tilde{u}_n(\overline{x}_n) \abs{\tilde{v}_n(\overline{x}_n)}^{\beta} \phi \left( \frac{\overline{x}_n -x_n}{r_n} \right) \d \overline{x}_n\\
    & = r_n^{-\frac{d-sp}{p}} \int_{\rd}  \abs{\tilde{u}_n(\overline{x}_n)}^{\al -2} \tilde{u}_n(\overline{x}_n) \abs{\tilde{v}_n(\overline{x}_n)}^{\beta} \phi \left( \frac{\overline{x}_n -x_n}{r_n} \right) \d \overline{x}_n \\
    & = \int_{\rd}  \abs{\tilde{u}_n(\overline{x}_n)}^{\al -2} \tilde{u}_n(\overline{x}_n) \abs{\tilde{v}_n(\overline{x}_n)}^{\beta} \phi_n(\overline{x}_n) \d \overline{x}_n.
\end{align*}
Further, applying Lemma \ref{convergence-integrals}-(i), we have 
\begin{align}\label{PSD-5}
    &\lim_{n \ra \infty} \int_{\rd} \abs{\hat{u}_n( x)}^{\al-2}\hat{u}_n(x) \abs{\hat{v}_n(x)}^{\beta} \phi(x) \dx = \int_{\rd} \abs{\hat{u}( x)}^{\al-2}\hat{u}(x) \abs{\hat{v}(x)}^{\beta} \phi(x) \dx, \no \\
    & \lim_{n \ra \infty} \int_{\rd} \abs{\hat{v}_n( x)}^{\be-2}\hat{v}_n(x) \abs{\hat{u}_n(x)}^{\al} \psi(x) \dx = \int_{\rd} \abs{\hat{v}( x)}^{\be-2}\hat{v}(x) \abs{\hat{u}(x)}^{\al} \psi(x) \dx.
\end{align}
Now taking the limit as $n \ra \infty$ in \eqref{PSD-4}, using \eqref{PSD_3},  \eqref{PSD-5} we obtain the claim.  

\noi \textbf{Step 4:} We set 
\begin{align*}
    &w_n(z) = \tilde{u}_n(z) - r_n^{-\frac{d-sp}{p}} \hat{u} \left( \frac{z-x_n}{r_n} \right), h_n(z) = \tilde{v}_n(z) - r_n^{-\frac{d-sp}{p}} \hat{v} \left( \frac{z-x_n}{r_n} \right), \\
    & \tilde{w}_n(z) = r_n^{\frac{d-sp}{p}} w_n(r_n z + x_n) \text{ and } \tilde{h}_n(z) = r_n^{\frac{d-sp}{p}} h_n(r_n z + x_n), \text{ for } z \in \rd.
\end{align*}
Note that $\norm{(w_n, h_n)}_{\W} = \norm{(\tilde{w}_n, \tilde{h}_n)}_{\W}$ and $\int_{\rd} \abs{w_n}^{\al} \abs{h_n}^{\be} \dx = \int_{\rd} \abs{\tilde{w}_n}^{\alpha} \abs{\tilde{h}_n}^{\be} \dx$. In this step, we show that $(w_n , h_n)$ is a (PS) sequence of $I_{0,0}$ at level $\eta - I_{f,g}(\tilde{u}, \tilde{v}) - I_{0,0}(\hat{u}, \hat{v})$. 
Observe that $\tilde{w}_n = \hat{u}_n - \hat{u}$ and $\tilde{h}_n = \hat{v}_n - \hat{v}$, and hence the norm invariance gives $\norm{(w_n, h_n)}_{\W} = \norm{(\hat{u}_n - \hat{u}, \hat{v}_n - \hat{v})}_{\W}$. Applying Lemma \ref{convergence2}, Lemma \ref{convergence-BL}, and Step 1 we calculate 
\begin{align*}
    &I_{0,0}(w_n, h_n) = \frac{1}{p} \norm{(\hat{u}_n- \hat{u}, \hat{v}_n - \hat{v})}_{\W}^p - \frac{1}{p^*_s} \int_{\rd}\abs{\hat{u}_n - \hat{u}}^{\alpha} \abs{\hat{u}_n - \hat{v}}^{\beta} \\
    & = \frac{1}{p} \left( \norm{\hat{u}_n}_{\wps}^p - \norm{\hat{u}}_{\wps}^p \right) + \frac{1}{p} \left(\norm{\hat{v}_n}_{\wps}^p - \norm{\hat{v}}_{\wps}^p \right) - \frac{1}{p^*_s}  \int_{\rd} \left( \abs{\hat{u}_n}^{\al} \abs{\hat{v}_n}^{\be} -  \abs{\hat{u}}^{\al} \abs{\hat{v}}^{\be}\right) + o_n(1) \\
    & = \frac{1}{p} \left( \norm{\tilde{u}_n}_{\wps}^p - \norm{\hat{u}}_{\wps}^p \right) + \frac{1}{p} \left( \norm{\tilde{v}_n}_{\wps}^p - \norm{\hat{v}}_{\wps}^p \right) - \frac{1}{p^*_s}  \int_{\rd} \left( \abs{\tilde{u}_n}^{\al} \abs{\tilde{v}_n}^{\be} -  \abs{\hat{u}}^{\al} \abs{\hat{v}}^{\be}\right) + o_n(1)\\
    & = I_{0,0}(\tilde{u}_n, \tilde{v}_n) - I_{0,0}(\hat{u}, \hat{v}) + o_n(1) = \eta - I_{f,g}(\tilde{u}, \tilde{v}) - I_{0,0}(\hat{u}, \hat{v}) + o_n(1).
\end{align*}
Next, we show $\prescript{}{\W'}{\langle} I_{0,0}'(w_n, h_n), (\phi, \psi) {\rangle}_{\W} \ra 0$ for every $\phi, \psi \in \wps$. Using the fact that $\{(w_n, h_n)\}$ is bounded on $\W$, we apply \eqref{HS2}, to obtain the following for some $C>0$:
\begin{align*}
    \left| \prescript{}{\W'}{\langle} I_{0,0}'(w_n, h_n), (\phi, \psi) {\rangle}_{\W} \right| \le C \norm{(\phi, \psi)}_{\W}, \text{ for every } (\phi, \psi) \in \W. 
\end{align*}
Hence using the density of $\cc (\rd)$ in $\wps$, and using the uniform boundedness principle, it is enough to show $\prescript{}{\W'}{\langle} I_{0,0}'(w_n, h_n), (\phi, \psi) {\rangle}_{\W} \ra 0$, for every $\phi, \psi \in \cc (\rd)$.
We take $(\phi,\psi)\in\cc (\rd)\times \cc (\rd)$ and then we consider
\begin{align*}
    \hat{\phi}_n(z) = r_n^{\frac{d-sp}{p}} \phi(r_n z + x_n), \text{ and } \hat{\psi}_n(z) = r_n^{\frac{d-sp}{p}} \psi(r_n z + x_n), \text{ for } z \in \rd.
\end{align*}
Since $\norm{(\hat{\phi}_n, \hat{\psi}_n)}_{\W} = \norm{(\phi, \psi)}_{\W}$, the sequence $\{ (\hat{\phi}_n, \hat{\psi}_n)\}$ is bounded on $\W$ and up to a subsequence $(\hat{\phi}_n, \hat{\psi}_n) \rightharpoonup (\hat{\phi}, \hat{\psi})$ in $\W$. Suppose $x_n \ra x_0$. Since $r_n \ra 0$ or $\infty$, in either case we get $(\hat{\phi}_n, \hat{\psi}_n)\rightharpoonup (0,0)$ in $\W$. Indeed, when $r_n\to 0$, $(\hat{\phi}_n,\hat{\psi}_n)\to (0,0)$ uniformly in $\R^{2d}$ and $\{ (\hat{\phi}_n, \hat{\psi}_n)\}$ is bounded in $\W$ enforce it has a weak limit (up to a subsequence) in $\W$ which must coincide with $(0,0)$. Therefore, $\hat{\phi}=0$ and $\hat{\psi} =0$ a.e. in $\rd$. Now, when $r_n\to\infty$, $\supp{\phi}$, and $\supp{\psi}$ contained in $B(0,R)$ for some $R>0$ will give $\supp{\hat{\phi}_n}$ and $\supp{\hat{\psi}_n}$ to concentrate at a point (maybe at infinity), which renders $ (\hat{\phi}_n, \hat{\psi}_n)\rightharpoonup (0,0)$ in $\W$. Further, we consider $|x_n| \ra \infty$ and $r_n \ra r_0 \in [0, \infty)$. When $r_0>0$, this is the case when the centre of the $L^{p^*_s}$ mass is escaping to infinity in a fixed radius and when $r_0= 0$, then fattening out $L^{p^*_s}$ mass and moving the center towards infinity are happening simultaneously. In all these cases, since $\|\hat{\phi}_n\|_{L^{p^*_s}} = \|\phi\|_{L^{p^{*}_s}}$ and $\|\hat{\psi}_n\|_{L^{p^*_s}} = \|\psi\|_{L^{p^{*}_s}}$, thus weak limit tends to 0 as $n\to\infty$. For $|x_n| \ra \infty$ and $r_n \ra \infty$, from the above argument, we similarly get the zero weak limit of $(\hat{\phi}_n, \hat{\psi}_n)$.  

Now we calculate 
\begin{align}\label{conv-1}
    \prescript{}{\W'}{\langle} I_{0,0}'(w_n, h_n), (\phi, \psi) {\rangle}_{\W} & =  \mathcal{A}(w_n, \phi) + \mathcal{A}(h_n, \psi) \nonumber \\
    &- \frac{\al}{p^*_s} \int_{\rd} \abs{w_n}^{\al -2} w_n \abs{h_n}^{\be} \phi - \frac{\be}{p^*_s} \int_{\rd} \abs{h_n}^{\be-2} h_n \abs{w_n}^{\al} \psi \nonumber \\
    & = \mathcal{A}(\tilde{w}_n, \hat{\phi}_n) + \mathcal{A}(\tilde{h}_n, \hat{\psi}_n) \nonumber \\
    &- \frac{\al}{p^*_s} \int_{\rd} \abs{\tilde{w}_n}^{\al -2} \tilde{w}_n \abs{\tilde{h}_n}^{\be} \hat{\phi}_n - \frac{\be}{p^*_s} \int_{\rd} \abs{\tilde{h}_n}^{\be-2} \tilde{h}_n \abs{\tilde{w}_n}^{\al} \hat{\psi}_n. 
\end{align}
By Lemma \ref{convergence-BL}-(ii), and using H\"{o}lder's inequality with the conjugate pair $(p, p')$ and further using $\norm{(\hat{\phi}_n, \hat{\psi}_n)}_{\W} = \norm{(\phi, \psi)}_{\W}$, we get 
\begin{align*}
 \mathcal{A}(\hat{u}_n - \hat{u}, \hat{\phi}_n) -  \mathcal{A}(\hat{u}_n, \hat{\phi}_n) + \mathcal{A}(\hat{u}, \hat{\phi}_n) = o_n(1), \; \mathcal{A}(\hat{v}_n - \hat{v}, \hat{\psi}_n) - \mathcal{A}(\hat{v}_n, \hat{\psi}_n) + \mathcal{A}(\hat{v}, \hat{\psi}_n) = o_n(1).
\end{align*}
Hence, using the change of variable, 
\begin{align*}
    \mathcal{A}(\tilde{w}_n, \hat{\phi}_n) = \mathcal{A}(\tilde{u}_n, \phi) - \mathcal{A}(\hat{u}, \hat{\phi}_n) + o_n(1), \,
    \mathcal{A}(\tilde{h}_n, \hat{\psi}_n) = \mathcal{A}(\tilde{v}_n, \psi) - \mathcal{A}(\hat{v}, \hat{\psi}_n) + o_n(1).
\end{align*}
Now we use the fact that $\tilde{u}_n, \tilde{v}_n, \hat{\phi}_n, \hat{\psi}_n$ weakly go to zero in $\wps$, and apply Lemma \ref{convergence-integrals} to get $\mathcal{A}(\tilde{u}_n, \phi) = o_n(1), \mathcal{A}(\tilde{v}_n, \psi) = o_n(1), \mathcal{A}(\hat{u}, \hat{\phi}_n)= o_n(1)$, and $\mathcal{A}(\hat{v}, \hat{\psi}_n) = o_n(1)$. Therefore, 
\begin{align*}
    \mathcal{A}(\tilde{w}_n, \hat{\phi}_n) = o_n(1) \text{ and } \mathcal{A}(\tilde{h}_n, \hat{\psi}_n) =o_n(1).
\end{align*}
Further, using a similar set of arguments (see Step 2), we also get 
\begin{align*}
   &\int_{\rd} \abs{\tilde{w}_n}^{\al -2} \tilde{w}_n \abs{\tilde{h}_n}^{\be} \hat{\phi}_n = \int_{\rd} \abs{\hat{u}_n - \hat{u}}^{\al -2} (\hat{u}_n - \hat{u}) \abs{\hat{v}_n - \hat{v}}^{\be} \hat{\phi}_n = o_n(1), \text{ and } \\
   &\int_{\rd} \abs{\tilde{h}_n}^{\be-2} \tilde{h}_n \abs{\tilde{w}_n}^{\al} \hat{\psi}_n =  o_n(1).
\end{align*}
Hence from \eqref{conv-1} we conclude $\lim_{n \ra \infty} \prescript{}{\W'}{\langle} I_{0,0}'(w_n, h_n), (\phi, \psi) {\rangle}_{\W} = 0,$ as required. Thus, by the density arguments, $(w_n , h_n)$ is a (PS) sequence for $I_{0,0}$. 

\noi \textbf{Step 5:} For $\hat{u},\, \hat{v},\, r_n,\, x_n$ as given in Step~2, we set $\tilde{u}_1 = \hat{u}, \tilde{v}_1 = \hat{v}, r_n^1 = r_n, x_n^1 =x_n$. Further, using a similar set of arguments as given in Step 2 and Step 3, there exist $r_n^2 \in (0, \infty), x_n^2 \in \rd$ such that if 
\begin{align*}
    \hat{w}_n(z) = (r_n^2)^{\frac{d-sp}{p}} w_n(r_n^2 z + x_n^2) \text{ and } \hat{h}_n(z) = (r_n^2)^{\frac{d-sp}{p}} h_n(r_n^2 z + x_n^2), \text{ for } z\in \rd,
\end{align*}
then $(\hat{w}_n, \hat{h}_n) \rightharpoonup (\hat{w}, \hat{h}) \neq (0,0)$ in $\W$, where $(\hat{w}, \hat{h})$ satisfies the homogeneous system \eqref{homogenous-system}. Moreover, the next sequence will satisfy the (PS) condition as level $\eta - I_{f,g}(\tilde{u}, \tilde{v}) - I_{0,0}(\hat{u}, \hat{v}) - I_{0,0}(\hat{w}, \hat{h})$. We set $\tilde{u}_2 = \hat{w}, \tilde{v}_2 = \hat{h}$.
Now in the same spirit of \cite{Tintarev} (Page 130, Theorem~3.3) and using Lemma \ref{Asy-Bub-inter}, we get 
\begin{align*}
    \mathcal{A}(C_{x^1_n,r^1_n} \tilde{u}_1,C_{x^2_n,r^2_n} \tilde{u}_2) = \mathcal{A}(C_{x^1_n,r^1_n} \hat u,C_{x^2_n,r^2_n} \hat w) = \mathcal{A}\left(\hat u, C_{\frac{x^2_n-x^1_n}{r_n^1}, \frac{r_n^2}{r_n^1}} \hat w \right) \ra 0\mbox{ as }n\to\infty. 
\end{align*}
Hence, in view of Proposition \ref{weak-bub}, we get
\bee\no
\bigg|\log\left(\frac{r^1_n}{r_n^2}\right)\bigg|+\bigg|\frac{x_n^1-x_n^2}{r_n^1}\bigg|\to\infty\quad \text{as }n\to\infty.
\eee 
Now, starting from a (PS) sequence $\{(\tilde u_n, \tilde v_n)\}$ for $I_{0,0}$ we have extracted further (PS) sequences at a level which is strictly lower than the level of $\{(\tilde u_n, \tilde v_n)\}$, and with a fixed amount of decrease in every step, since ($(\hat u,\hat v)$ is a solution to \eqref{homogenous-system}), 
\begin{align*}
    J_{0,0}(\hat u,\hat v) =\frac{s}{d}\|(\hat u,\hat v)\|_{\W}^p \text{ and } S_{\alpha,\beta}\leq \|(\hat u,\hat v)\|_{\W}^{\frac{sp^2}{d}},
\end{align*}
which in turn implies $J_{0,0}(\hat u,\hat v)\geq \frac{s}{d}S_{\alpha,\beta}^{\frac{d}{sp}}.$ Similarly, $J_{0,0}(\hat w,\hat h)\geq \frac{s}{d}S_{\alpha,\beta}^{\frac{d}{sp}}.$ 
On the other hand, since we have $\sup_n\|(\tilde u_n, \tilde v_n)\|_{\W}$ is finite, there exists $k \in \N$ such that this process terminates after the $k$ number of steps and the last (PS) sequence strongly converges to $0$. Moreover, 
\bee\no
\bigg|\log\left(\frac{r^i_n}{r_n^j}\right)\bigg|+\bigg|\frac{x_n^i-x_n^j}{r_n^i}\bigg|\to\infty\quad \text{for } i\neq j,\, \, 1\leq i,\, j\leq k.
\eee
This completes the proof.
\qed

\section{Existence of a positive solution for non-homogeneous system}

For $f,\, g$ non-negative nontrivial functionals in $(\wps)'$ with $\ker (f) = \ker (g)$, we introduce the following functional:
\begin{align}\label{20-24-5'}
J_{f,g}(u,v)\coloneqq \frac{1}{p}\|(u,v)\|^p_{\mathcal{W}}
-\frac{1}{p^*_s}\int_{\rd} u_+^{\alpha}v_+^{\beta} \dx & - \prescript{}{(\wps)'}{\langle}f,u{\rangle}_{\wps} \no \\
&-\prescript{}{(\wps)'}{\langle}g,v{\rangle}_{\wps}, \; \forall\, (u,v) \in \W. 
\end{align}
Observe that $J_{f,g} \in \C^1\big(\mathcal{W}\big)$. Suppose, $(\bar u, \bar v) \in \W$ is a nontrivial critical point of $J_{f,g}$, then it satisfies the following system of equations weakly:
\begin{equation}\tag{$\mathcal S\geq 0$}\label{hom-sys-pos}
\left\{\begin{aligned}
&(-\Delta_p)^s u = \frac{\al}{p_{s}^{*}(m)}u_{+}^{\alpha-1}v_{+}^{\beta}+f(x)\;\;\text{in}\;\rd,\\
&(-\Delta_p)^s v = \frac{\beta}{p^{*}_{s}(m)}v_{+}^{\beta-1}u_{+}^{\alpha}+g(x)\;\;\text{in}\;\rd,
\end{aligned}
\right.
\end{equation}
i.e., the following identity holds
\begin{align}\label{system_eqn"}
   &\iint\limits_{\rd\times\rd}\frac{|\bar{u}(x)-\bar{u}(y)|^{p-2}(\bar{u}(x)-\bar{u}(y))(\phi(x)-\phi(y))}{|x-y|^{d+sp}} \dxy \no \\ 
   &+\iint\limits_{\rd\times\rd}\frac{|\bar{v}(x) - \bar{v}(y)|^{p-2}(\bar{v}(x)-\bar{v}(y))(\psi(x)-\psi(y))}{|x-y|^{d+sp}} \dxy \no \\ 
   &=\frac{\al}{p^*_s} \int_{\rd} \bar{u}_{+}^{\alpha-1}\bar{v}_{+}^{\beta} \phi \dx + \frac{\be}{p^*_s} \int_{\rd} \bar{v}_{+}^{\beta-1}\bar{u}_{+}^{\alpha} \psi \dx + \prescript{}{(\wps)'}{\langle}f,\phi{\rangle}_{\wps}+ \prescript{}{(\wps)'}{\langle}g,\psi{\rangle}_{\wps}, 
\end{align}
for every $(\phi, \psi) \in \W$. Taking $(\bar{u}_{-},\bar{v}_{-})$ as a test function in \eqref{system_eqn"} we get
\begin{align*}
   &-\iint\limits_{\rd\times\rd}\frac{|\bar{u}(x)-\bar{u}(y)|^{p-2}\left\{(\bar{u}_{+}(x)\bar{u}_{-}(y)+\bar{u}_{+}(y)\bar{u}_{-}(x)) + (\bar{u}_{-}(x)-\bar{u}_{-}(y))^2\right\}}{|x-y|^{d+sp}} \dxy\\
   &\qquad\qquad- \iint\limits_{\rd\times\rd}\frac{|\bar{v}(x)-\bar{v}(y)|^{p-2}\left\{(\bar{v}_{+}(x)\bar{v}_{-}(y)+\bar{v}_{+}(y)\bar{v}_{-}(x)) + (\bar{v}_{-}(x)-\bar{v}_{-}(y))^2\right\}}{|x-y|^{d+sp}} \dxy\\
   &=\prescript{}{(\wps)'}{\langle}f,\bar{u}_{-}{\rangle}_{\wps}+ \prescript{}{(\wps)'}{\langle}g,\bar{v}_{-}{\rangle}_{\wps} \geq 0.
\end{align*}
This in turn implies $\prescript{}{(\wps)'}{\langle}f,\bar{u}_{-}{\rangle}_{\wps} = 0 = \prescript{}{(\wps)'}{\langle}g,\bar{v}_{-}{\rangle}_{\wps}$. Taking $(-\bar{u}_{-}, 0)$ as a test function in \eqref{system_eqn"}, we obtain
\begin{align*}
   \norm{\bar{u}_-}_{\wps}^p = - \prescript{}{(\wps)'}{\langle}f,\bar{u}_{-}{\rangle}_{\wps} = 0,
\end{align*}
which implies $\bar u_-=0$ a.e. in $\rd$. Choosing $(0, -\bar{v}_{-})$, we also get $\bar v_-=0$ a.e. in $\rd$. Therefore, $(\bar u, \bar v)$ is a nonnegative solution of~\eqref{system_eqn}. Next, we assert that $(\bar u, \bar v)\neq (0,0)$ implies $\bar u\neq 0$ and $\bar v\neq 0$.
Suppose not, that is, assume, for instance, $\bar u\neq 0$ but $\bar v=0$. Then taking $(\bar u, 0)$ and $(0, \bar u)$ as tests function in \eqref{system_eqn"} we get 
$$ \norm{\bar u}_{\wps}^p=\prescript{}{(\wps)'}{\langle}f,\bar u{\rangle}_{\wps}, \text{ and } \prescript{}{(\wps)'}{\langle}g, \bar u{\rangle}_{\wps}=0.$$
Hence we get $\bar u\in \ker (f)=\ker (g)$ and $\norm{\bar u}_{\wps}^p=0$, which in turn implies $\bar u=0$. This contradicts the fact that $(\bar u, \bar v)\neq (0,0)$. Similarly, if $\bar u=0$, then $\bar v=0$ too. Thus, the assertion holds. Further, taking $(\phi, 0)$ with $\phi \ge 0$ as a test function in \eqref{system_eqn"} we get 

\begin{align*}
    \iint\limits_{\rd\times\rd}\frac{|\bar{u}(x)-\bar{u}(y)|^{p-2}(\bar{u}(x)-\bar{u}(y))(\phi(x)-\phi(y))}{|x-y|^{d+sp}} \dxy & = \frac{\al}{p^*_s} \int_{\rd} \bar{u}^{\alpha-1}\bar{v}^{\beta} \phi \dx \\
    &+ \prescript{}{(\wps)'}{\langle}f,\phi{\rangle}_{\wps} \ge 0.
\end{align*}
Now, applying the strong maximum principle for the $(s,p)-$super harmonic function (see \cite[Theorem 1.2]{DQ}), we obtain $\overline{u} > 0$ a.e. in $\rd$. Similarly, $\overline{v} > 0$ a.e. in $\rd$. 
Next, suppose $\bar u \equiv \bar v$. First, we assume $f\equiv g$ but $\alpha\neq \beta$. Then, taking the test function $(\bar u,-\bar u)$ in \eqref{system_eqn"} yields
$$\frac{1}{p^*_s}(\alpha-\beta)\int_{\rd} \bar{u}^{\alpha+\beta} \,{\rm d}x=0.$$
This is impossible since $\bar u$ is positive in $\rd$. In the remaining case, we assume $\alpha=\beta$ but $f\not\equiv g$ and $\text{ker}(f)=\text{ker}(g)$. Then taking the test function $(\phi,-\phi)$, where $\phi\in \wps$, we obtain
$$ \prescript{}{(\wps)'}{\langle}f-g,\phi{\rangle}_{\wps}=0.$$
This in turn implies $f\equiv g$ as $\phi\in \wps$ is arbitrary. 

In view of the above discussion, we have the following remark.

\begin{remark}\label{rmk-1}
    Let $f, g \in (\wps)'$ where $f,\, g$ is nonnegative nontrivial functionals and $\ker (f) = \ker (g)$. If $(\overline{u}, \overline{v}) \in \W$ is a  critical point of $J_{f, g}$, then $\overline{u}, \overline{v} > 0$ a.e. in $\rd$. Moreover, $\bar{u} \neq \bar{v}$, whenever $f\equiv g$ but $\alpha \neq \beta$; or $\alpha = \beta$ but $f \not \equiv g$. Hence, if $(u,v)\in\W$ is a critical point of $J_{f,g}$ then $(u,v)$ is a solution of \eqref{system_eqn}.
\end{remark}

To prove the existence of  critical point for $J_{f,g}$, we consider the following sets 
\begin{align*}
    &\Upsilon_1\coloneqq\left\{(u,v)\in\mathcal{W}\,:\, (u,v)=(0,0) \,\,\text{or}\,\, \Psi(u,v)>0\right\}, \\
    &\Upsilon\coloneqq\left\{(u,v)\in\mathcal{W} \setminus \{(0,0)\}\,:\, \Psi(u,v)=0\right\}, \text{ and } \\
    &\Upsilon_2\coloneqq\left\{(u,v)\in\mathcal{W}: \Psi(u,v)<0\right\},
\end{align*}
via a functional $\Psi:\mathcal{W}\to\mathbb{R}$ defined by
\begin{align*}
    \Psi(u,v)\coloneqq\norm{(u,v)}_{\mathcal{W}}^p-\left(\frac{p_s^*-1}{p-1}\right)\int_{\rd}|u|^{\al}|v|^{\beta}\dx, \; \forall \, (u,v) \in \W.
\end{align*}
Note that $\Psi \in \C^1\big(\mathcal{W}\big)$.
\begin{remark}\label{omega-1}
(i)  For $t>0$, notice that 
\begin{align}\label{t-1}
    \Psi(tu,tv)=0 \iff t= \left(\frac{(p-1)\|(u,v)\|_{\W}^p}{(p_s^*-1)\int_{\rd}|u|^{\alpha}|v|^{\beta}\,{\rm d}x}\right)^{\frac{1}{p_s^*-p}}.
\end{align}
If  $(u,v) \in \Upsilon_1\setminus\{(0,0)\}$, then using the fact that $\|(u,v)\|_{\W}^p>\left(\frac{p_s^*-1}{p-1}\right)\int_{\rd}|u|^{\alpha}|v|^{\beta}\,{\rm d}x$, there is $t_1>1$ (given in \eqref{t-1}) such that $(t_1u,t_1v)\in \Upsilon$. Similarly, if $(u,v)\in \Upsilon_2\setminus\{(0,0)\}$, there is $t_2<1$ (given in \eqref{t-1}) such that $(t_2u,t_2v)\in \Upsilon$.

(ii) By observing that
$$\Psi(t u,t v)=(t^p-t^{p_s^*})\|(u,v)\|_{\W}^p, \; \forall \,(u,v)\in \Upsilon,$$ 
we get $(t u,t v)\in \Upsilon_1$ for all $t\in(0,1)$ and  and $(t u,t v)\in \Upsilon_2$ for all $t>1.$

\end{remark}

Consider the following quantities 
\begin{align*}
    c_0\coloneqq\inf_{\Upsilon_1}J_{f,g}(u,v), \text{ and } c_1 \coloneqq\inf_{\Upsilon}J_{f,g}(u,v).
\end{align*}
\begin{lemma}\label{strict}
    If 
    \begin{align*}
        \max\{\|f\|_{(\wps)'}, \|g\|_{(\wps)'}\} \le \frac{1}{2} \left( \frac{1}{p} - \frac{p-1}{p^*_s(p^*_s-1)} \right) \left(  \frac{p-1}{p^*_s-1} S_{\al, \be}^{\frac{p^*_s}{p}} \right)^{\frac{p-1}{p^*_s-p}},
    \end{align*}
    where $S_{\al, \be}$ is the best constant of \eqref{HS2}, then $c_0<c_1$. 
\end{lemma}

\begin{proof}
   For brevity, we call the R.H.S. by $\tilde{\var}\equiv \tilde{\var}(d,p,s)$. If $(u,v) \in \Upsilon$, then   
   \begin{align*}
       \norm{(u,v)}_{\mathcal{W}}^p =\left(\frac{p_s^*-1}{p-1}\right)\int_{\rd}|u|^{\al}|v|^{\beta}\dx.
   \end{align*}
   Let $S_{\al, \be}$ be the best constant of \eqref{HS2}. Then for every $(u,v) \in \Upsilon$, the above identity gives 
   \begin{align*}
       \norm{(u,v)}_{\W}^{p^*_s-p} \ge \frac{p-1}{p^*_s-1} S_{\al, \be}^{\frac{p^*_s}{p}} \Longrightarrow \norm{(u,v)}_{\W} \ge \left(  \frac{p-1}{p^*_s-1} S_{\al, \be}^{\frac{p^*_s}{p}} \right)^{\frac{1}{p^*_s-p}} =: \delta \left(\equiv\delta(d,p,s)\right).
   \end{align*}
   Moreover, for every $(u,v) \in \Upsilon$, 
   \begin{align*}
       J_{f,g}(u,v) & = \left( \frac{1}{p} - \frac{p-1}{p^*_s(p^*_s-1)} \right) \norm{(u,v)}_{\W}^p - \prescript{}{(\wps)'}{\langle}f,u{\rangle}_{\wps} - \prescript{}{(\wps)'}{\langle}g,v{\rangle}_{\wps} \\
       & \ge \norm{(u,v)}_{\W} \left( \left( \frac{1}{p} - \frac{p-1}{p^*_s(p^*_s-1)} \right) \norm{(u,v)}_{\W}^{p-1} - \norm{f}_{(\wps)'} - \norm{g}_{(\wps)'} \right) \\
       & \ge \delta  \left( \left( \frac{1}{p} - \frac{p-1}{p^*_s(p^*_s-1)} \right) \delta^{p-1} - \norm{f}_{(\wps)'} - \norm{g}_{(\wps)'} \right) \\
       &= \delta \left( 2\tilde{\var}  - \norm{f}_{(\wps)'} - \norm{g}_{(\wps)'} \right) \ge 0,
   \end{align*}
   which implies $c_1\ge 0$. Next, from 
   \begin{align*}
& J_{f,g}(tu,tv) = \frac{t^p}{p}\|(u,v)\|^p_{\mathcal{W}}
-\frac{t^{p^*_s}}{p^*_s}\int_{\rd} u_+^{\alpha}v_+^{\beta} \dx - t\prescript{}{(\wps)'}{\langle}f,u{\rangle}_{\wps} -t\prescript{}{(\wps)'}{\langle}g,v{\rangle}_{\wps}, \text{ and } \\
& \Psi(tu,tv) = t^p \norm{(u,v)}_{\mathcal{W}}^p- t^{p^*_s} \left(\frac{p_s^*-1}{p-1}\right)\int_{\rd}|u|^{\al}|v|^{\beta}\dx, 
   \end{align*}
   we note that $J_{f,g}(tu,tv) <0$ and $\Psi(tu,tv)>0$ for every $t <<1$. Therefore, $c_0 \le J_{f,g}(tu,tv) < 0 \le c_1$, as required.
\end{proof}

\begin{proposition}\label{fir-cri-pt}
Let $s \in(0, \frac{1}{p'})$. Put
\begin{align*}
c_0=\inf_{\Upsilon_1}J_{f,g}(u,v).    
\end{align*}
If $\max\{\|f\|_{(\wps)'}, \|g\|_{(\wps)'}\} \le \tilde{\var}$, then $J_{f,g}$ has a critical point $(u_0,v_0)\in \overline{\Upsilon}_1$ with $J_{f,g}(u_0,v_0)=c_0<0$. 
\end{proposition}

\begin{proof} We decompose the proof into a few steps.
		
\noi \textbf{Step~1:} In this step, we see $c_0>-\infty$. 
Since $J_{f,g}(u,v)\geq I_{f,g}(u,v)$, it is enough to show that $I_{f,g}$ is bounded from below. From the definition of $\Upsilon_1$, it immediately follows that for all $(u,v) \in \Upsilon_1$,
\begin{align}\label{31-7-3}
    I_{f,g}(u,v)\geq \left( \frac{1}{p}-\frac{p-1}{p^*_s(p^*_s-1)} \right) \|(u,v)\|_{\W}^p- \left( \|f\|_{(\wps)'}+\|g\|_{(\wps)'} \right) \|(u, v)\|_{\W}.
\end{align}
As R.H.S. is $p^{\text{th}}$ power function in $\|(u,v)\|_{\W}$ (for $p>1$), $I_{f,g}$ is bounded from below. Hence, Step~$1$ follows. 

\noi \textbf{Step~2:} In this step, we show that there exists a bounded nonnegative (PS)-sequence $\{(u_n,v_n)\} \subset \Upsilon_1$ for $J_{f,g}$ at level $c_0$. Let $\{(u_n,v_n)\}\subset \overline{\Upsilon}_1$ such that $J_{f,g}(u_n,v_n)\to c_0$. Applying Lemma \ref{strict} and Ekeland's variational principle, from $\{(u_n,v_n)\} $ we can extract a (PS) sequence in $\Upsilon_1$ for $J_{f,g}$ at level $c_0$. We again call it by $\{(u_n,v_n)\}$. Moreover, as $J_{f,g}(u,v)\geq I_{f,g}(u,v)$, from \eqref{31-7-3} it follows that $\{(u_n,v_n)\} $ is a bounded sequence. Therefore,  up to a subsequence $(u_n,v_n)\rightharpoonup (u_0,v_0)$ in $\W$ and $(u_n,v_n)\to (u_0,v_0)$ a.e. in $\rd$. In particular, $ (u_n)_+\to (u_0)_+,\;(v_n)_+\to (v_0)_+$ and $ (u_n)_-\to (u_0)_-, \;(v_n)_-\to (v_0)_-$ a.e. in $\rd$. Moreover, as $f,\;g$ are nonnegative functionals, a straight-forward computation yields
\begin{align*}
    o_n(1)&=\prescript{}{\W'}\langle J'_{f,g}(u_n,v_n), \left((u_n)_-,(v_n)_-\right)\rangle_{\W}\\
		&=\mathcal{A}(u_n,(u_n)_-)+\mathcal{A}(v_n,(v_n)_-)-\prescript{}{(\wps)'}\langle f, (u_n)_-\rangle_{\wps}-\prescript{}{(\wps)'}\langle g,(v_n)_-\rangle_{\wps}\\
		&\leq-\|\left((u_n)_-,(v_n)_-\right)\|_{\W}^p,\no
\end{align*}
where the final inequality follows using the following elementary inequality (see \cite[Lemma A.2]{BrPa}):
\begin{align*}
    \abs{a - b}^{p-2}(a-b)(a_- - b_-) \le -\abs{a_- - b_-}^{p}, \text{ for } a,b \in \R.
\end{align*}
Therefore, $\big((u_n)_-,(v_n)_-\big)\to (0,0)$ in $\W$, which implies up to a subsequence $(u_n)_-\to 0$ and $(v_n)_-\to 0$ a.e. in $\rd$ and thus $(u_0)_-=0$  and $(v_0)_-=0$ a.e. in $\rd$.  Consequently, without loss of generality, we can assume that $\{(u_n,v_n)\}$ is a nonnegative sequence. This completes the proof of Step~2.

\noi \textbf{Step 3:}  In this step we show that $(u_n,v_n)\to (u_0,v_0)$ in $\W$ and $(u_0,v_0)\in \bar \Upsilon$. Applying Theorem \ref{PS-decomposition}, we get
\bee\label{J5}
(u_n,v_n) \to  (u_0,v_0) +\sum_{j=1}^{k} (C_{x_n^j, r_n^j}\tilde{u}_j, C_{x_n^j, r_n^j}\tilde{v}_j), \;\text{ in }\W,\text{ as }n\to\infty,
\eee
with $J_{f,g}'(u_0,v_0) =0$, $(\tilde u_j,\tilde v_j)$  is a nonnegative solution of \eqref{homogenous-system} ($(u_n, v_n)$ is (PS) sequence for $J_{f,g}$ implies $(\tilde u_j,\tilde v_j)$  a solution of \eqref{hom-sys-pos} with $f=g=0$ and by Remark~\ref{rmk-1}, $(\tilde u_j,\tilde v_j)$  is a nonnegative solution of \eqref{homogenous-system}), and  $\{x_n^j\}_n \subset \rd, \{r_n^j\}_n \subset \R^+$ are some appropriate sequences such that either $x_n^j\to x^j$ and $r_n^j \ra 0, \text{or } \infty$, or $|x_n^j|\to\infty$, as $n\to\infty$. To prove Step~3, we need to show that $k=0$. Arguing by contradiction, suppose that $k\neq 0$ in \eqref{J5}. Therefore, for $1\leq j\leq k$ we have
		\begin{align}\label{12-4-3}
		\Psi\bigg((C_{x_n^j, r_n^j}\tilde{u}_j, C_{x_n^j, r_n^j}\tilde{v}_j)\bigg)&=\|\left(\tilde u_j,\tilde v_j\right)\|^p_{\W}-\left(\frac{p_s^*-1}{p-1}\right)\int_{\rd}|\tilde u_j|^{\alpha}|\tilde{v}_j|^{\beta}\dx\no\\
        &=\left(\frac{p-p_s^*}{p-1}\right) \|\left(\tilde u_j,\tilde v_j\right)\|_{\W}^p<0.
		\end{align}			
		From Theorem \ref{PS-decomposition}, we also have
		$$c_0=\lim_{n\to\infty} J_{f,g}(u_n,v_n)=  J_{f,g}(u_0,v_0)+\sum_{j=1}^{k}J_{0,0}(\tilde u_j,\tilde v_j).$$
Since $(\tilde u_j,\tilde v_j)$ is a solution to \eqref{homogenous-system}, we see that $J_{0,0}(\tilde u_j,\tilde v_j)\geq \frac{s}{d}S_{\alpha,\beta}^{\frac{d}{sp}}.$
Consequently, $ J_{f,g}(u_0,v_0)<c_0$. Therefore, $(u_0,v_0)\not\in \Upsilon_1$ (as $\inf_{\Upsilon_1} J_{f,g}(u,v)=:c_0$) and
\bee\label{12-4-4} \Psi(u_0,v_0)\leq 0.\eee 	
Next, we evaluate $\Psi\bigg((u_0,v_0) +\sum_{j=1}^{k}(C_{x_n^j, r_n^j}\tilde{u}_j, C_{x_n^j, r_n^j}\tilde{v}_j)\bigg)$. We observe that $(u_n,v_n)\in \Upsilon$ implies $\Psi(u_n,v_n)\geq 0$. Combining this with the uniform continuity of $\Psi$ and \eqref{J5} yields
\bee\label{J8}
0\leq \liminf_{n\rightarrow\infty}\Psi(u_n,v_n)=\liminf_{n\rightarrow\infty} \Psi\bigg((u_0,v_0) +\sum_{j=1}^{k} (C_{x_n^j, r_n^j}\tilde{u}_j, C_{x_n^j, r_n^j}\tilde{v}_j) \bigg).
\eee
Note that from Step~2, we already have $u_0,v_0\geq 0$ and $\tilde{u}_j, \tilde{v}_j \ge 0$ for all $j$. Therefore, 
\begin{align}\label{12-4-1}
&\Psi\bigg((u_0,v_0) +\sum_{j=1}^{k} (C_{x_n^j, r_n^j}\tilde{u}_j, C_{x_n^j, r_n^j}\tilde{v}_j) \bigg)\no\\
&=\bigg{\|}(u_0,v_0)+\sum_{j=1}^{k}(C_{x_n^j, r_n^j}\tilde{u}_j, C_{x_n^j, r_n^j}\tilde{v}_j)\bigg{\|}^p_{\W}\no\\
&\qquad- \left(\frac{p_s^*-1}{p-1}\right)\int_{\rd}\bigg|u_0+\sum_{j=1}^{k} C_{x_n^j, r_n^j}\tilde{u}_j   \bigg|^{\alpha}\bigg|v_0+\sum_{j=1}^{k} C_{x_n^j, r_n^j}\tilde{v}_j \bigg|^{\beta}\dx\no\\
&\underbrace{\leq}_{\text{Claim}} \Psi(u_0,v_0)+ \sum_{j=1}^{k}\Psi \left( C_{x_n^j, r_n^j}\tilde{u}_j, C_{x_n^j, r_n^j}\tilde{v}_j \right)+ o_n(1).
\end{align}
We have
\bee\no
\int_{\Rn}\bigg|u_0+\sum_{j=1}^{m}C_{x_n^j, r_n^j}\tilde{u}_j   \bigg|^{\al}\bigg|v_0+\sum_{j=1}^{m}C_{x_n^j, r_n^j}\tilde{v}_j\bigg|^{\beta}\dx\geq \bigg(\int_{\Rn} u_0^\al v_0^\beta{\rm d}x+\sum_{j=1}^{m}\int_{\Rn} \left| C_{x_n^j, r_n^j}\tilde{u}_j \right|^\al \left| C_{x_n^j, r_n^j}\tilde{v}_j \right|^\beta {\rm d}x \bigg).
\eee
To prove the claim in \eqref{12-4-1}, we need to show the following: 
\begin{align}\label{split-1}
   \bigg{\|}(u_0,v_0)+\sum_{j=1}^{k}(C_{x_n^j, r_n^j}\tilde{u}_j, C_{x_n^j, r_n^j}\tilde{v}_j)\bigg{\|}^p_{\W}- \|(u_0,v_0)\|_{\W}^p- \sum_{j=1}^{k} \bigg{\|}(C_{x_n^j, r_n^j}\tilde{u}_j, C_{x_n^j, r_n^j}\tilde{v}_j) \bigg{\|}^p_{\W} = o_n(1). 
\end{align}
First, we show that 
\begin{align}\label{split-2}
  \bigg{\|} u_0 + \sum_{j=1}^{k}C_{x_n^j, r_n^j}\tilde{u}_j \bigg{\|}_{\wps}^p  - \norm{u_0}_{\wps}^p - \sum_{j=1}^{k} \bigg{\|} C_{x_n^j, r_n^j}\tilde{u}_j \bigg{\|}_{\wps}^p = o_n(1). 
\end{align}
To prove \eqref{split-2} we use the following inequality (see \cite[Lemma~1.2 of Technical Lemmas]{Bahri}):
\bee\label{fie}
\Bigg{|}\bigg{|}\sum_{j=0}^{k}a_j\bigg{|}^{p-1}\sum_{j=0}^{k}a_j - \sum_{j=0}^{k}\abs{a_j}^{p-1}{a_j}\Bigg{|} \lesssim_{p} \sum_{0\leq i\neq j\leq k}\abs{a_j}^{p-1}\abs{a_i},
\eee
with
\begin{align*}
   a_j= C_{x_n^j, r_n^j}\tilde{u}_j(x)-C_{x_n^j, r_n^j}\tilde{u}_j(y),\text{ where, (as notation) }C_{x_n^0, r_n^0}\tilde{u}_0\coloneqq u_0. 
\end{align*}
Then using \eqref{fie} we have 
\begin{align}\label{ele-1}
   &\Bigg{|}\sum_{j=0}^{k} \left( C_{x_n^j, r_n^j}\tilde{u}_j(x)-C_{x_n^j, r_n^j}\tilde{u}_j(y) \right) \Bigg{|}^p - \sum_{j=0}^{k} \bigg{|}C_{x_n^j, r_n^j}\tilde{u}_j(x)-C_{x_n^j, r_n^j}\tilde{u}_j(y)\bigg{|}^p \no \\
    &= \left| \sum_{j=0}^{k} a_j\right|^p - \sum_{j=0}^{k} \abs{a_j}^p \le \Bigg{|}\bigg{|}\sum_{j=0}^{k}a_j\bigg{|}^{p-1}\sum_{j=0}^{k}a_j - \sum_{j=0}^{k}\abs{a_j}^{p-1}{a_j}\Bigg{|} \lesssim_{p} \sum_{0\leq i\neq j\leq k}\abs{a_j}^{p-1}\abs{a_i} \no \\
    &= C(p) \sum_{0\leq i\neq j\leq k} \left| C_{x_n^j, r_n^j}\tilde{u}_j(x)-C_{x_n^j, r_n^j}\tilde{u}_j(y) \right|^{p-1}\left|C_{x_n^i, r_n^i}\tilde{u}_i(x)-C_{x_n^i, r_n^i}\tilde{u}_i(y)\right|,
\end{align}
for some constant $C(p)>0$. The first inequality in \eqref{ele-1} follows from 
\begin{align*}
    &\abs{l} - \sum_{j=0}^{k} \abs{a_j}^p \le \abs{l} - \abs{m} \le \abs{\abs{l} - \abs{m}} \le \abs{l-m}, \text{ with } \\
    &l= \bigg{|}\sum_{j=0}^{k}a_j\bigg{|}^{p-1}\sum_{j=0}^{k}a_j, \text{ and } m=\sum_{j=0}^{k}\abs{a_j}^{p-1}{a_j}.
\end{align*}
Hence using \eqref{ele-1},
\begin{align*}
    &\text{L.H.S of \eqref{split-2}} \\
    &\lesssim_p \sum_{0\leq i\neq j\leq k} \;\iint\limits_{\rd\times\rd} \frac{\left| C_{x_n^j, r_n^j}\tilde{u}_j(x)-C_{x_n^j, r_n^j}\tilde{u}_j(y) \right|^{p-1}\left|C_{x_n^i, r_n^i}\tilde{u}_i(x)-C_{x_n^i, r_n^i}\tilde{u}_i(y)\right|}{|x-y|^{d+sp}} \dxy.
\end{align*}
In view of Theorem \ref{PS-decomposition}, 
\begin{align*}
    \left| \log \left( \frac{r_n^i}{r_n^j} \right) \right| + \left| \frac{x_n^i - x_n^j}{r_n^i}  \right| \rightarrow \infty, \text{ for } i \neq j, 1 \le i, j \le k.
\end{align*}
Therefore, applying Lemma \ref{weak-bub-converse-III} leads to 
\begin{align*}
    \sum_{0\leq i\neq j\leq k} \;\iint\limits_{\rd\times\rd} \frac{\left|C_{x_n^j, r_n^j}\tilde{u}_j(x)-C_{x_n^j, r_n^j}\tilde{u}_j(y)\right|^{p-1}\left|C_{x_n^i, r_n^i}\tilde{u}_i(x)-C_{x_n^i, r_n^i}\tilde{u}_i(y)\right|}{|x-y|^{d+sp}} \dxy = o_n(1).
\end{align*}
Thus, \eqref{split-2} holds true. Similarly, 
\begin{align*}
    \bigg{\|} v_0 + \sum_{j=1}^{k}C_{x_n^j, r_n^j}\tilde{v}_j \bigg{\|}_{\wps}^p  - \norm{v_0}_{\wps}^p - \sum_{j=1}^{k} \bigg{\|} C_{x_n^j, r_n^j}\tilde{v}_j \bigg{\|}_{\wps}^p = o_n(1). 
\end{align*}
As a consequence, we obtain \eqref{split-1}. Combining \eqref{12-4-3} and \eqref{12-4-4} with \eqref{12-4-1} we get a contradiction to \eqref{J8} for large enough $n \in \N$. Therefore, $k=0$ in \eqref{J5}, which implies $(u_n,v_n)\to (u_0,v_0)$ in $\W$ as $n\to\infty$. Consequently, $\Psi(u_n,v_n)\to \Psi(u_0,v_0)$ as $n\to\infty$, which in turn implies $(u_0,v_0)\in \bar \Upsilon$. Thus Step~3 follows.

\noi {\bf Step 4:} From the previous steps we see that $J_{f,g}(u_0,v_0)=c_0$ and $J_{f,g}'(u_0,v_0)=0$. Therefore, $(u_0,v_0)$ is a weak solution to \eqref{hom-sys-pos}. Finally, in view of the proof of Lemma \ref{strict}, $c_0<0.$
\end{proof}

\noi {\bf Proof of Theorem~\ref{th:mul}:} From Proposition~\ref{fir-cri-pt} combining with Remark~\ref{rmk-1} the proof of Theorem~\ref{th:mul} follows. \qed

\begin{remark}
    Initially, our main goal was to prove the existence of at least two positive solutions of \eqref{system_eqn} in the spirit of \cite{MoSoMiPa}. Indeed, provided that $\norm{f}_{(\wps)'}, \,\norm{g}_{(\wps)'}$ are sufficiently small, we obtain $c_0<c_1$ where
    \bee\no
    c_1\coloneqq\inf_{\Upsilon}J_{f,g}(u,v).
    \eee
    But unlike the case $p=2$, it is not immediate for us to split $\|(u_1+u_2,v_1+v_2)\|_{\W}^{p}$-norm in such a way (see proof of \cite[Proposition 4.2]{MoSoMiPa}) that we could create a Mountain Pass critical level higher than $c_1$ and lower than the first level of breaking down (PS) condition (where the (PS) condition of $I_{f,g}$ holds), which in turn gives rise to a second positive solution.
\end{remark}

\noi \textbf{Acknowledgments:}
N. Biswas acknowledges the support of the Science and Engineering Research Board, Government of India, for National Postdoctoral Fellowship (PDF/2023/000038). S. Chakraborty acknowledges the support of the Tata Institute of Fundamental Research, Bengaluru for the Institute Postdoctoral Fellowship. The authors express their gratitude to Prof. Lorenzo Brasco for the valuable discussions that contributed to the writing of this manuscript. We are also grateful to the anonymous referees for their insightful comments and valuable suggestions, which significantly improved the clarity and presentation of the paper.

\noi \textbf{Declaration:} 
The authors declare that no data was generated or analyzed during the current study. As such, there is no data availability associated with this manuscript. Furthermore, the authors declare that no conflict of interest concerning the authorship of this manuscript.

\bibliographystyle{abbrvnat}

\end{document}